\documentclass[a4paper,12pt]{scrartcl}
\usepackage[english]{babel}
\usepackage[T1]{fontenc}
\usepackage[utf8]{inputenc}
\DeclareUnicodeCharacter{00A0}{ }
\DeclareUnicodeCharacter{014C}{\=O}
\usepackage{mathptmx}
\usepackage[scaled=.92]{helvet}
\usepackage{courier}
\usepackage{version}
\usepackage{paralist}
\usepackage[nodetail]{optional}
\usepackage{amsmath,amsfonts,amssymb}

\newcommand{\RN}[1]{\uppercase\expandafter{\romannumeral #1\relax}}

\newcommand{\N}{\mathbb{N}}

\newcommand{\R}{\mathbb{R}}
\newcommand{\C}{\mathbb{C}}
\newcommand{\RE}{\operatorname{Re}}
\newcommand{\IM}{\operatorname{Im}}
\newcommand{\Hom}{\operatorname{Hom}}
\newcommand{\End}{\operatorname{End}}
\newcommand{\divergence}{\operatorname{div}}
\newcommand{\eps}{\varepsilon}
\newcommand{\hh}{\hslash}
\newcommand{\intd}{\mathop{\operatorname{d}}}
\newcommand{\supp}{\operatorname{supp}}

\newcommand{\Id}{\operatorname{Id}}

\newcommand{\volg}{\mu_g}
\newcommand{\tr}{\operatorname{tr}}
\newcommand{\spec}{\operatorname{spec}}
\newcommand{\outgoing}{{\mathrm{out}}}
\newcommand{\incoming}{{\mathrm{in}}}
\newcommand{\suchthat}{\mathrel{;}}
\newcommand{\bigoh}{\mathcal{O}}
\newcommand{\kernel}{\operatorname{ker}}

\newcommand{\image}{\operatorname{im}}
\newcommand{\inner}{\,}
\newcommand{\Ell}{\ensuremath{\mathcal E}}
\newcommand{\Mix}{\ensuremath{\mathcal M}}
\newcommand{\Hyp}{\ensuremath{\mathcal H}}
\newcommand{\Char}{\operatorname{Char}}
\newcommand{\nozerosection}{\setminus 0}
\newcommand{\Cinfty}{\ensuremath{C^{\infty}}}
\newcommand{\Ccinfty}{\ensuremath{\Cinfty_c}} 
\newcommand{\Dprime}{{\ensuremath{\mathcal{D}^\prime}}} 
\newcommand{\op}{\operatorname{Op}} 
\newcommand{\WF}{\operatorname{WF}}
\newcommand{\bWF}{\mathbin{^b \operatorname{WF}}}
\newcommand{\bTstar}{\mathbin{^b T^*}} 
\newcommand{\Sphg}{S_{\mathrm{phg}}} 
\newcommand{\Psiphg}{\Psi_{\mathrm{phg}}} 
\newcommand{\Psitang}{\Psi_{\mathrm{tang}}} 
\newcommand{\hnabla}{\mathcal{H}} 
\newcommand{\vnabla}{\mathcal{V}} 
\newcommand{\Lie}{\mathcal{L}} 
\newcommand{\PT}{\tau} 
\newcommand{\PTE}{\PT^E} 
\newcommand{\from}{\gets} 
\newcommand{\EE}{E\times E}
\newcommand{\DB}[1]{|T#1|} 
\newcommand{\hDB}[1]{|T#1|^{1/2}} 

\usepackage{amsthm}
\newtheorem{theorem}{Theorem}[section]
\newtheorem{proposition}[theorem]{Proposition}
\newtheorem{corollary}[theorem]{Corollary}
\newtheorem{lemma}[theorem]{Lemma}
\theoremstyle{definition}

\newtheorem{remark}[theorem]{Remark}

\usepackage[pdftex]{color}
\usepackage[pdftex]{graphicx}
\usepackage[pdftex,
	colorlinks=true,
	linkcolor=gruen,
        urlcolor=blue,
        citecolor=blue,
	bookmarksopen=false,
	pdfpagemode=UseNone]{hyperref}
\definecolor{gruen}{rgb}{0,0.5,0}
\hypersetup{pdfauthor={S.~Hansen},pdftitle={Propagation of Polarization in Transmission Problems}}

\title{Propagation of Polarization \\ in Transmission Problems}
\author{Sönke~Hansen%
\date{November 24, 2021}
\footnote{Universität Paderborn, Institut für Mathematik, Warburger Str.~100, 33098 Paderborn}}
\makeindex
\begin{document}

\maketitle
\begin{abstract}
For geometric systems of real principal type, we define a subprincipal symbol
and derive a transport equation for polarizations which, in the scalar case,
is a well-known equation of Duistermaat and Hörmander.
We apply the transport equation to propagation of polarization in transmission
problems of elastodynamics, to interior bulk waves as well as to free (Rayleigh) surface waves.
Using spectral factorizations of matrix polynomials having real spectrum,
we establish reflection and refraction laws of polarizations at the boundary and at interior interfaces.
The results are not limited to isotropic elasticity.
{}
\end{abstract}

\section{Introduction} 
Polarization of an elastic or an electromagnetic wave
refers to the vector-valued amplitude in the high-frequency regime.
Extending methods of geometrical optics, microlocal analysis applies to prove
well-posedness of propagation of polarizations along rays broken by
reflection and transmission (refraction) at boundaries and at interior interfaces.
Such results are fundamental to geometric methods for solving the inverse problem
of recovering interfaces and elasticities.
Only recently, in \cite{StefUhlVasy21transm}, has microlocal well-posedness of
the transmission problem of isotropic elastodynamics been proven.
The purpose of this paper is to prove, also for anisotropic elastic media,
propagation of polarizations along broken bicharacteristics.
Making assumptions about sources, we exclude glancing rays from consideration.

The elasticity operator $L$ only depends on the metric tensor and the stiffness $4$-tensor of the elastic material.
Away from sources, elastic waves satisfy the homogeneous wave equation $Lu=\rho D_t^2 u$,
where $D_t=i^{-1}\partial/\partial t$, $i$ the imaginary unit, and $\rho>0$ is the material density.
Our basic assumption is that the operator of elastodynamics, $P=L-\rho D_t^2$, is of real principal type.
This covers isotropic and generic elastic media.

For systems of real principal type, $Pu=0$, the propagation of polarization wavefront sets in the interior
was studied in \cite{Dencker82polarrpt}, and at boundaries in \cite{Gerard85polar}.
Suppose $u$ is a Lagrangian distribution section in a bundle $E$, $\Lambda$ the associated Lagrangian manifold.
Then one expects finer results.
In particular, the polarization, that is the principal symbol $a$ of $u$,
should satisfy a transport equation.
In the scalar case, \cite{DuisHorm72FIOtwo} derived the transport equation
$i^{-1}\Lie_H a + p_s a = 0$, where $\Lie_H$ is
the Lie derivative in direction of the Hamilton field $H$ and $p_s$ the subprincipal symbol of $P$.
For general systems with non-scalar principal symbol,
no definition of subprincipal symbol and no transport equation seems to be known.
Suppose the manifold and the bundle $E$ are endowed with connections $\nabla$ and $\nabla^E$.
In Theorem~\ref{theorem-amp-transport-eqn} we show that the polarization $a$ satisfies a transport equation
\begin{equation}
\label{eq-transport-polariz}
i^{-1}D_H a + p_s a = 0 \quad\text{on $\Lambda$.}
\end{equation}
The first order differential operator $D_H$ is assembled via Leibniz' formula
from $\nabla^E$ pulled back to $\Lambda$ and the Lie derivative $\Lie_H$ on half-densities.
The geometric pseudo-differential calculus of \cite{Sharaf05geosymb} is basic to
the definition of the subprincipal symbol $p_s$ and to the proof of \eqref{eq-transport-polariz}.
For outgoing solutions of $Pu=f$ with Lagrangian sources $f$,
we apply the Lagrangian intersection calculus of \cite{MelroseUhlmann79intersection}
to derive initial conditions for the polarization.

An elastic body is modelled over a three-dimensional Riemannian manifold $M$ with boundary $\partial M$.
The free surface problem asks for elastic waves $u$ having zero boundary traction,
$Tu=0$ at $\partial M\times \R_t$.
Layers having different elastic properties are separated by hypersurfaces $N$ in the interior of $M$.
In classical formulation, the transmission problem requires displacements $u$
and tractions to be continuous across $N\times \R_t$.
In this paper we study the transmission problem of elastodynamics,
looking for outgoing waves generated by Lagrangian sources located in the interior,
on the boundary, or on the interfaces.

It is standard to reduce the analysis of reflection of wavefront sets at a boundary
to the construction of parametrices for Dirichlet data.
Finding such parametrices is possible after a suitable
spectral decomposition on the principal symbol level, \cite{Taylor75reflection}.
We establish, microlocally near a non-glancing boundary region, factorizations
\[ L -\rho D_t^2 = (D_r -Q^\sharp)A_0 (D_r - Q). \]
Here $r\geq 0$ is the distance to the boundary, and the first order pseudodifferential
operators $Q$ and $Q^\sharp$ are tangential, that is, they commute with multiplication by $r$.
The spectrum of the principal symbol $q$ of $Q$ is contained in the closure of the complex upper half-plane.
Following \cite{Taylor75reflection}, the equation $(D_r-Q)u=0$ is solved for arbitrary initial data at $r=0$.
On the principal symbol level, the factorization of the elastodynamics operator corresponds to,
and the operator factorization is derived from,
spectral factorizations of self-adjoint quadratic polynomials in $s\in\C$:
\begin{equation}
\label{eq-specfact-of-eladyn-symbol}
\ell(\eta+s\nu)-\rho\tau^2\Id =(s-q^\sharp(\eta,\tau))a_0(s-q(\eta,\tau)),
\end{equation}
$\eta,\nu=\intd r(y)\in T_y^*\partial M$.
The polynomial takes its values in $\End(E)$.
Here the Hilbert space $E$ is the complexification of the tangent space of $M$ at $y\in\partial M$.
The spectra of $q$ and $q^\sharp$ are disjoint.
There holds $E=E_c\oplus E_r$, where $E_c$ and $E_r$ are the sums of generalized eigenspaces of $q$
associated to eigenvalues $s$ satisfying $\IM s>0$ and $s\in\R$, respectively.
In case $E_r=0$, the right root $q$ is unique.
On the other hand, if $E_r\neq 0$, then, because of a sign characteristic of real eigenvalues,
there are two distinct right roots $q$, $q_\outgoing$ and $q_\incoming$.
It is a consequence of the sign characteristic combined with the real principal type property that,
along the bicharacteristics of $D_r-Q_\outgoing$ (resp.\ of $D_r-Q_\incoming$)
which issue from the boundary into the interior, time $t$ increases (resp.\ decreases).
Using Dirichlet parametrices,
\[ (D_r- Q_{\outgoing/\incoming}) U_{\outgoing/\incoming}\equiv 0, \quad U_{\outgoing/\incoming}|_{r=0}\equiv\Id, \]
we introduce the DN (Dirichlet to Neumann) operators $Z_{\outgoing/\incoming}= TU_{\outgoing/\incoming}$.
If $Z_\outgoing$ is elliptic, reflection is given by $Z_{\outgoing}^{-1}Z_\incoming$ in the
sense that this operator maps, microlocally, incoming to outgoing Dirichlet data.
Similarly, reflection and transmission at interior interfaces is readily
resolved provided $Z^+_\outgoing+Z^-_\outgoing$ is elliptic;
the superscripts refer to sides of the interface.

The principal symbol $z=z_{\outgoing}$ of $Z_\outgoing$ satisfies the useful inequality
\begin{equation}
\label{eq-energy-flow-ineq}
-\tau \IM(zv|v) >0 \quad\text{if $v\not\in E_c$, and $\geq 0$ for every $v\in E$.}
\end{equation}
The inequality reflects, on the principal symbol level, outgoing energy flow.
The proof is based on the spectral factorization \eqref{eq-specfact-of-eladyn-symbol}.
The inequality \eqref{eq-energy-flow-ineq} implies $\kernel z\subset E_c$ and,
at the interface, $\kernel(z^++z^-)\subset E_c^+\cap E_c^-$.
It follows that $Z_\outgoing$ and $Z^+_\outgoing+Z^-_\outgoing$ are elliptic operators
in the hyperbolic region, see Proposition~\ref{prop-z-ell-at-Hyp}.
We emphasize that, at the interface, the hyperbolic region $\Hyp$ defined in this paper
contains any region which is hyperbolic with respect to at least one side of the interface.
For isotropic media,
we show by explicit symbol computations that the operators are elliptic also in mixed regions,
Proposition~\ref{prop-z-ell-at-Mix}.

We include an analysis of the propagation of polarizations of free (Rayleigh) surface waves
which may occur over the elliptic boundary region.
Rayleigh wave propagation was recognized as a propagation of singularities phenomenon
by \cite{Taylor79rayleighwaves}.
The theory of free surface waves in anisotropic elastic media, \cite{LotheBarnett85surfwaveimped},
implies that restrictions to the boundary of Rayleigh waves satisfy a real principal type system,
\cite{Nakamura91rayleighpulses},
thereby making these waves accessible to refined microlocal analysis.

In the final Section~\ref{section-propag-polar} we state and prove results on
propagation of polarization originating from Lagrangian sources.
We refrain from stating the results in terms of FIOs.

Pseudo-differential calculus with differential geometric structure is
important in the present work because it offers an invariant leading symbol
which includes the subprincipal level; in fact, it even has an invariant full symbol.
Formulas for the leading symbols of compositions and adjoints contain vertical
and horizontal derivatives of principal symbols.
See Section~\ref{section-horizontal-deriv} for details on these derivatives.
Section~\ref{section-geom-symbol-calc}
is a self-contained presentation of the geometric pseudo-differential calculus
of \cite{Sharaf05geosymb} to the extent needed in our paper.

There are two appendices.
Appendix~\ref{sect-elastic-symmetries} recalls elastic symmetries.
Furthermore, it is used to show that generically the elastodynamic operator
of tranversely isotropic perturbations of isotropic media is of real principal type.
In Appendix~\ref{sect-spectral-factorization} we redevelop, simplifying standard proofs,
the theory of spectral factorization of self-adjoint quadratic matrix polynomials,
and we deduce some corollaries.

\section{Linear elasticity}
\label{sect-linear-elast}
We present the equations of linear elasticity in differential geometric language,
and we apply the variational solution method and regularity theory to the
transmission problem of elastodynamics.

Let $M$ be a smooth ($\Cinfty$) manifold endowed with a smooth Riemannian metric tensor $g$
and the associated Levi-Civita connection $\nabla$.
Denote $T^{(p,q)}M$ the bundle of $p$-contravariant and $q$-covariant tensors.
So, $TM=T^{(1,0)}M$ is the tangent and $T^*M=T^{(0,1)}M$ the cotangent bundle.
Furthermore, $\End(TM)=T^{(1,1)}M$.
The Lie derivative of $g$ along a given vector field $u\in\Cinfty(M;TM)$
is a symmetric tensor field $\Lie_u g\in \Cinfty(M;T^{(0,2)}M)$.
Identifying the tensor bundles $T^{(0,2)}M$ and $T^{(1,1)}M$ using $g$,
$\Lie_u g$ agrees with the symmetrization of $\nabla u\in\Cinfty(M;\End(TM))$.

For an elastic material modeled on $(M,g)$, infinitesimal displacements $u$ are vector fields.
The tensor field $\eps(u) =\Lie_u g$ is the strain induced by the displacement $u$.
Elastic materials are distinguished by their stiffness (or elasticity) tensors which are $4$-tensor
fields $C$ which map strains to stresses $\sigma=C\eps$.
In coordinate notation $\sigma^{ij}=C^{ijkm}\eps_{km}$.
(Here we use the summation convention.)
We regard the stiffness tensor as an endomorphism of $T^{(1,1)}M$
as well as an homomorphism from $T^{(0,2)}M$ into its dual $T^{(2,0)}$.
A stiffness tensor field $C$ is assumed to satisfy the following properties:
$C$ is a symmetric homomorphism which annihilates antisymmetric tensors, maps into symmetric tensors,
and is positive definite when restricted to the subspace of symmetric tensors.
The last condition is known in elasticity as the strong convexity condition.
When $\dim M=3$, stiffness tensors are classified by their $SO(3)$ symmetry subgroups.
Isotropy and transverse isotropy are two of eight distinct symmetry classes; see Appendix~\ref{sect-elastic-symmetries}.
We say that $(M,g,C)$ is an elastic body if $M$ is an oriented, compact, connected
manifold with smooth boundary, $\partial M\neq\emptyset$.

Let $(M,g,C)$ be an elastic body.
Suppose $C\in\Cinfty$.
Denote by $x\inner \bar y$ the inner product induced by $g$ on tangent and tensor spaces and their complexifications;
the bar marks the conjugate linear slot.
Let $\intd V_M$ and $\intd V_{\partial M}$ denote the Riemannian volume elements on $M$ and on its boundary.
The elasticity operator is the second order differential operator $L$ given by
the following identity of Green's type:
\begin{equation}
\label{eq-Green-id-elasticity}
\int_M C\eps(u)\inner\overline{\eps(v)}\;\intd V_M
  =\int_M L u \inner\bar{v}\;\intd V_M + \int_{\partial M} T u \inner\bar{v}\;\intd V_{\partial M}.
\end{equation}
The traction operator $T$ is a first order differential operator.
Using local coordinates, we denote by ${}_{;j}$ covariant derivatives with respect to the $j$-th coordinate,
and we lower indices to turn $u$ into a covector.
Then $\eps_{km}=(u_{k;m}+u_{m;k})/2$.
By \eqref{eq-Green-id-elasticity} and the symmetries of $C$, and using the divergence theorem, we have
\begin{equation*}
(Lu)^i = -(C^{ijkm} u_{k;m})_{;j},
\quad
(Tu)^i = C^{ijkm}u_{k;m} \nu_j.
\end{equation*}
Here $\nu$ is the exterior unit conormal at $\partial M$;
$\nu=-\intd r|_{\partial M}$ with $0\leq r$ the distance to $\partial M$.
The principal symbols of $L$ and $T$ at a covector $\xi\in T^*M$ are given by the matrices
$\ell(\xi)=(C^{ijkm}\xi_j\xi_m)$ and $(\sqrt{-1}C^{ijkm}\nu_j\xi_m)$, respectively.
Observe that $L$ is uniformly elliptic.

By Korn's inequality and the positive definiteness of $C$ the sesquilinear form
\[ E(u,v)= \int_M C\eps(u)\inner\overline{\eps(v)}\;\intd V_M =\int_M C\nabla u\inner\overline{\nabla v}\;\intd V_M  \]
is coercive on the Sobolev space $H^1=H^1(M;\C\otimes TM)$.
Following the variational method, the weak definition of the elasticity operator $L:H^1\to{H^1}^*$ is given by
\begin{equation}
\label{eq-weak-elastic-traction-problem}
E(u,v)= (Lu|v), \quad u,v\in H^1.
\end{equation}
The duality bracket is induced via the $L^2$ inner product: $(w|v)=\int_M w\inner\bar v\;\intd V_M$.
For $\lambda > 0$ sufficiently large, $L+\lambda$ is an isomorphism from $H^1$ onto its conjugate dual ${H^1}^*$.
Since $C$ is real symmetric, $L$ restricts to a self-adjoint operator on $L^2$ with domain
\[ D(L)=\{u\suchthat Lu\in L^2\} = \{u\in H^2\suchthat Tu|_{\partial M} =0\}.\]
The last equality follows from $H^2$-regularity combined
with Green's identity~\eqref{eq-Green-id-elasticity} for $u\in H^2$ and $v\in H^1$.
The self-adjoint operator $L$ is the elasticity operator of the free surface problem.
Variational regularity theory applies to give higher order Sobolev regularity:
If $u\in D(L)$ and $Lu\in H^k$, then $u\in H^{k+2}$.

More generally, we consider composite elastic bodies with components separated by smooth hypersurfaces.
To be precise, suppose that $C$ is smooth except for possible jump discontinuities across
a smooth orientiable hypersurface $N$ contained in the interior $M^\circ$ of $M$.
The hypersurface $N$ is not connected in general.
The components of $M':=M^\circ\setminus N$ represent the interiors of subbodies of $M$, layers for example.
The boundaries of subbodies are submanifolds of $\partial M\cup N$ equipped with their orientations as boundaries.
In the following, Sobolev spaces are of integer order and consist of sections of the complexified tangent bundle;
we omit the bundles from the notation.
The definition \eqref{eq-weak-elastic-traction-problem} of the elasticity operator $L$ is still applicable.
There holds $H^2$-regularity:
If $Lu\in L^2(M)$, then $u\in H^1(M)\cap H^2(M')$.
For $u\in H^2(M')$, define the traction jump at $N$ as $[Tu]=(Tu)|_{N_+}+(Tu)|_{N_-}$;
the exterior unit normals of the sides $N_\pm$ of $N$ satisfy $\nu_+=-\nu_-$.
Applying \eqref{eq-Green-id-elasticity} to each component of $M'$ and summing over the components, we get
\[
\int_M C\eps(u)\inner\overline{\eps(v)}\;\intd V_M
  =\int_M L u \inner\bar{v}\;\intd V_M + \int_{\partial M} T u \inner\bar{v}\;\intd V_{\partial M}
    + \int_{N} [Tu] \inner\bar{v}\;\intd V_{N}.
\]
The restriction of $L$ to the domain
\[ D(L_T) = \{u\suchthat Lu\in L^2(M)\} = \{u\in H^1(M)\cap H^2(M')\suchthat Tu|_{\partial M}=0,\; [Tu]=0\}. \]
is a non-negative self-adjoint operator $L_T$, the elasticity operator corresponding to
homogeneous traction and transmission conditions.
For $u\in H^1(M')$ the jump $[u]=u|_{N_+}-u|_{N_-}$ of $u$ across $N$ vanishes iff $u\in H^1(M)$.
So, $D(L_T)$ consists of all $u\in H^2(M')$ which satisfy the zero traction boundary condition,
$Tu|_{\partial M}=0$, and the homogeneous transmission conditions at the interior interfaces, $[u]=[Tu]=0$ at $N$.
To determine the sign of $[u]$, fix an orientation of $N$.
Higher regularity holds: If $u\in D(L)$ and $Lu\in H^k(M')$, then $u\in H^{k+2}(M')$.
Moreover, the inverse $(L_T+1)^{-1}$ exists and maps $H^k(M')$ onto $H^{k+2}(M')\cap D(L_T)$.

Now we consider elastic waves $u=u(x,t)$ over $M\times \R_t$.
Put $D_t=-\sqrt{-1}\partial/\partial_t$.
Given sources $f$, $h$, and $h_j$ over $M'\times\R$, $\partial M\times\R$, and $N\times\R$, respectively,
we look for solutions of the following problem of elastodynamics:
\begin{equation}
\label{eq-eladyn-source-problem}
Lu-\rho D_t^2 u=f, \quad Tu=h, \quad [u]=h_0,\quad [Tu]=h_1.
\end{equation}
The material density $\rho>0$ is a smooth function on $M$ except for a jump discontinuity at $N$.
We are interested in outgoing solutions of \eqref{eq-eladyn-source-problem},
that is $\supp(u)\subset M\times[t_0,\infty[$ should hold for some time $t_0$.
In the main part of the paper we construct, with the help of microlocal parametrices,
and under assumptions on the wavefront sets of the sources,
approximate solutions which satisfy \eqref{eq-eladyn-source-problem} up to $\Cinfty$ errors.
To correct for these errors we need to solve \eqref{eq-eladyn-source-problem} with $\Cinfty$ data $f,h,h_j$.
If $r$ is a defining function of (part of) $\partial M\cup N$, then the commutator $[T,r]$ is non-singular.
Therefore, for $h$ and $h_j$ smooth, we can find $w=w(x,t)$ which is smooth in $M\times \R$
except for a jump discontinuity at $N\times \R$ such that $Tw=h$ holds at $\partial M$,
and $[w]=h_0$ and $[Tw]=h_1$ hold at $N\times \R$.
Therefore, we may assume $h=0=h_j$ and $f$ smooth, except for a jump at $N$.

We use the variational form of the evolution equation \eqref{eq-eladyn-source-problem},
and we apply standard theory, \cite{Evans98pde,Wloka87pde}, mutatis mutandis. 
Over a time intervall $t_0<t<t_1$ and given data $f\in L^2(t_0,t_1;L^2(M))$, we consider
\begin{equation}
\label{eq-eladyn-weak-form}
(u''(t)|v) +E(u(t),v)=(f(t)|v) \quad\text{for each $v\in H^1(M)$.}
\end{equation}
A prime denotes a time derivative.
In \eqref{eq-eladyn-weak-form} we absorbed $\rho$ into the $L^2$ inner product
of the Gelfand triplet $H^1(M)\hookrightarrow L^2(M)\hookrightarrow {H^1}^*(M)$
by redefining $(w|v)=\int_M w\inner\bar{v}\,\rho\intd V_M$, that is, we replaced the volume density by a mass density.
Furthermore, we replace $L$, $L_T$, and $f$ by the quotients $\rho^{-1}L$, $\rho^{-1}L_T$, and $f/\rho$.
There exists a unique solution $u(t)$ of \eqref{eq-eladyn-weak-form},
\[ u\in L^2(t_0,t_1;H^1(M)), \quad u'\in L^2(t_0,t_1;L^2(M)), \quad u''\in L^2(t_0,t_1;{H^1}^*(M)), \]
such that $u(t_0)=0$ and $u'(t_0)=0$ hold.
Initial values are defined because $u\in C([t_0,t_1];L^2(M))$ and $u'\in C([t_0,t_1];{H^1}^*(M))$.
Suppose $f\in H^k(t_0,t_1;L^2(M))$ satisfies the compatibility conditions between
initial and boundary values; for example, $f^{(j)}(t_0)=0$ holds for $0\leq j<k$.
Then we have additional regularity with respect to $t$:
\[ u\in H^k(t_0,t_1;H^1(M)), \quad u^{(k+1)}\in L^2(t_0,t_1;L^2(M)), \quad u^{(k+2)}\in L^2(t_0,t_1;{H^1}^*(M)). \]
Regularity with respect to the spatial variables follows using the differential equation
\[ u(t) = (L_T+1)^{-1} \big(f(t)- u''(t) + u(t)\big). \]
We can now combine the regularity results to conclude that, given a positive integer $m$,
there exists $k\geq m$ such that the following holds:
If $f\in H^k(t_0,t_1;H^k(M'))$ satisfies $f^{(j)}(t_0)=0$ for $j<k$, then the unique
solution $u$ of \eqref{eq-eladyn-weak-form} and $u(t_0)=u'(t_0)=0$ satisfies
\[ u\in H^m(t_0,t_1;H^m(M')\cap D(L_T)). \]
Thus, if $f$ is smooth in $]-\infty,t_1]$ and supported in $[t_0,t_1]$, then so is $u$.

\section{Covariant and horizontal derivatives}
\label{section-horizontal-deriv}
We need some facts, mostly standard \cite{KobaNomizu63I,Tu17diffgeom},
from the differential geometry of vector bundles.

Let $X$ be a $\Cinfty$ manifold endowed with a symmetric linear connection $\nabla$.
Assume $X$ Hausdorff, second countable, and without boundary.
The exponential map $\exp$ maps an open neighborhood of the zero-section of the tangent bundle $TX\to X$ into $X$.
We write $\exp_x v=\exp(v)$ for $v\in T_xX$.
A normal neighborhood centered at $x\in X$ is the
diffeomorphic image $\exp_x(V)$ of a star-shaped zero-neighborhood $V\subset T_xX$.
An open set which is a normal neighborhood centered at each of its points is called normal convex.
Any two points in a normal convex set $O$ are the endpoints of a unique, up to parametrization, geodesic in $O$.
If $O$ is a normal convex open set, then $(x,v)\mapsto (x,y)$, $y=\exp_x v$, maps
an open subset of $TX$ diffeomorphically onto $O\times O$.
The topology of $X$ has a basis consisting of normal convex sets.
See \cite[Ch.~\RN{1} \S 6]{Helgason01diffgeom}.

A local coordinate system $(x^j)$ on $X$ defines local frame fields $(\partial_j)$ of $TX$
and $(\intd x^j)$ of the cotangent bundle $\pi:T^*X\to X$.
Here we introduced $\partial_j$ as an abbreviation of $\partial/\partial x^j$.
The covariant derivative is given by
$\nabla\partial_j =\Gamma_{kj}^i\partial_i\intd x^k$ where $\Gamma_{kj}^i$ are the Christoffel symbols.
(We use the summation convention of summing over equal indices in opposite position.)
In the cotangent bundle, canonical coordinates
$(x^j,\xi_j)$ are defined by $\xi=\xi_j\intd x^j\in T_x^*X$.
We abbreviate $\partial/\partial\xi_j$ as $\partial^j$.

The Hessian of a function $\varphi\in\Cinfty(X)$ is given by
\[
\nabla^2\varphi=\nabla\intd\varphi
   = (\partial_k\partial_m \varphi-\Gamma_{km}^j \partial_j\varphi)\intd x^k \intd x^m.
\]
The symmetry of $\nabla$ implies that the Hessian is a symmetric $2$-tensor field.

Let $a\in\Cinfty(T^* X)$.
Using canonical coordinates, the vertical derivative $\vnabla a$ and the horizontal derivative $\hnabla a$
are defined as follows:
$\vnabla a(x,\xi)= \partial^j a(x,\xi)\partial_j$, and
\begin{equation}
\label{eq-horiz-a-scalar}
\hnabla a(x,\xi)=  \big(\partial_k a(x,\xi) + \Gamma^j_{km}(x)\xi_j \partial^m a(x,\xi)\big)\intd x^k.
\end{equation}
Invariantly, $\vnabla a$ and $\hnabla a$ are sections of the complexifications
of pullback bundles $\pi^* TX$ and $\pi^* T^*X$, respectively.
The vertical derivative $\vnabla a(x,\xi)\in\Hom(T_x^*X,\C)=T_xX\otimes\C$
is the ordinary derivative in fiber direction.
The horizontal derivative satisfies, and is uniquely determined by, the identity
\[
\intd(a\circ\intd\varphi)=(\hnabla a)\circ\intd\varphi + \nabla^2\varphi \cdot (\vnabla a\circ\intd\varphi),
\]
which holds for $\varphi\in\Cinfty(X;\R)$.
A centered dot denotes contraction of $T^*X\otimes TX$ factors to scalars by taking the trace.

Let $(E,\nabla^E)$ be a real or complex vector bundle $p^E:E\to X$ with a linear connection.
If $(e_j)$ is a local frame field of $E$, then the covariant derivative operator $\nabla^E$ 
is given by a matrix $(\omega^k_j)$ of connection one-forms:
\begin{equation}
\label{eq-E-frame-connection}
\nabla^E s = \intd s^j\otimes e_j + s^k\omega^{j}_k \otimes e_j
  \quad\text{if $s(x)=s^j(x) e_j(x)$,}
\end{equation}
or $\nabla^{E}_V s = (V s^j) e_j + s^k \omega^j_k(V) e_j$
for every vector field $V$.
If $f:Y\to X$ is smooth, then the pullback bundle $f^*E\to Y$ together with
a pullback connection $\nabla^{f^*E}$ are defined.
It is convenient to regard sections of $f^*E$ as sections along $f$,
that is as smooth maps $s:Y\to E$ which satisfy $p^E\circ s=f$.
Using the notation of \eqref{eq-E-frame-connection}, $(e_j\circ f)$ is a local frame field of $f^*E$,
and $f^* \omega^j_k$ are the connection one-forms of $\nabla^{f^*E}$:
\[
\nabla^{f^*E}_V s = (V s^j) e_j\circ f + s^k \omega^j_k(f_* V) e_j\circ f
  \quad\text{if $s(y)=s^j(y) e_j(f(y))$.}
\]
For example, suppose $f:\R\times X\to X$ is the projection along $\R$, $f(t,x)=x$,
and $W=a\partial_t+V$, $V$ tangent to $X$.
Then $\nabla^{f^*E}_W =a\partial_t + \nabla^E_V$.
The dual bundle $E^*$ and the density bundles $|E|^\alpha$ are equipped with
connections naturally induced from $\nabla^E$.
Suppose $(F,\nabla^F)$ is another vector bundle with a linear connection.
Then $E\otimes F$ and $\Hom(E,F)=E^*\otimes F$ are equipped 
with the linear connections induced by $\nabla^E$ and $\nabla^F$.
A connection on $\Hom(E,F)\to X$ is determined by the Leibniz rule:
\[ \nabla_V^F Au= (\nabla_V^{\Hom(E,F)}A)u+A\nabla_V^E u \]
holds if $u$ is a section of $E$, $A$ a section of $\Hom(E,F)$, and $V$ a vector field.

Suppose $E$ is equipped with a Hermitian metric $(\cdot|\cdot)_E$
which is compatible with $\nabla^E$, i.e.,
$V (u|v)_E=(\nabla^E_V u|v)_E+(u|\nabla^E_V v)_E$
holds for real vector fields $V$ and sections $u$ and $v$ of $E$.
The connection $\nabla^E$ is said to be being metric.
Hermitian metrics are linear in the first slot and antilinear in the second.

Suppose $X$ oriented and (pseudo-)Riemannian with metric tensor $g$,
$\nabla$ the Levi-Civita connection.
Then there is a canonical volume form $0<\volg\in\Cinfty(X;\DB{X})$
which is parallel, that is, $\nabla\volg=0$ holds.
The divergence $\divergence V$ of a vector field $V$ is defined by $\Lie_V\volg=(\divergence V)\volg$,
where $\Lie_V$ denotes the Lie derivative defined by $V$.
There holds
\begin{equation}
\label{eq-int-LieV}
\int_X (Vc+c\divergence V)\volg=\int_X \Lie_V (c\volg)=0 \quad\text{if $c\in\Ccinfty(X)$.}
\end{equation}
Furthermore, $\divergence V=\tr\nabla V$.
Consider the scalar product $\int_X (u|v)_E\volg$ of $u,v\in\Ccinfty(X;E)$.
It follows from \eqref{eq-int-LieV} that $-i\nabla^E_V -i\divergence V$ is the formal adjoint of $-i\nabla^E_V$.

Let $c:I\subset\R\to X$ be a smooth curve.
By definition, a section $s$ along $c$ of $E$ is parallel iff, on the interval $I$, $\nabla^{c^*E}s=0$ holds.
These equations are a homogeneous linear system of ordinary differential equations for $s$, i.e.,
\begin{equation}
\label{eq-ode-par-transp}
\dot s^j + \omega^j_k(c(t))\dot c(t) s^k =0
\quad\text{where $s(t)=s^j(t)e_j(c(t))$.}
\end{equation}
Dots denote derivatives with respect to $t$.
The fundamental matrices 
$\PTE_{c,t_0,t}:E_{c(t_0)}\to E_{c(t)}$ assign $s(t)$ to $s(t_0)$ when $\nabla^{c^*E}s=0$.
The linear map $\PTE_{c,t_0,t}$ is called the parallel transport in $E$ along $c$ from $c(t_0)$ to $c(t)$.
We abbreviate $\PTE_{c,t_0,t}$ as $\PTE_{c}$ if the endpoint parameters $t_0$ and $t$ are clear from the context.
If $y=\exp_x v$ belongs to a normal neighborhood centered at $x$,
then we let $\PTE_{x\from y}:E_y\to E_x$ denote the parallel transport along a
geodesic $c$ from $y$ to $x$, $c(t)=\exp_x(v-tv)$ for $0\leq t\leq 1$.
A reparametrization of a curve from $y$ to $x$ does not affect the parallel transport map $E_y\to E_x$.
Parallel transport along piecewise smooth curves is defined in an obvious way.
Covariant derivatives are recovered from the associated parallel transport:
\[ \big(\nabla^E_{\dot c(0)}s\circ c\big)(0) = \frac{\intd}{\intd t}\PTE_{c,t,0}s(c(t))|_{t=0}. \]

A loop based at $x$ is a curve with initial and end point equal to $x$.
For infinitesimal loops based at $x$, holonomy theory relates parallel transport to curvature.
We need a special case where the curvature does not contribute.
\begin{lemma}
\label{lemma-trivial-holonomy}
Let $c:[0,\delta[\to X$ be a smooth curve.
Set $x=c(0)$.
For $0<\eps<\delta$, denote by $\lambda_\eps$ the loop at $x$ which consists of $c|_{[0,\eps]}$
followed by the geodesic $x\from c(\eps)$.
Then $\PTE_{\lambda_\eps}=\Id_{E_x}+\bigoh(\eps^2)$ as $\eps\to 0$.
\end{lemma}
\begin{proof}
Assume that $c$ stays in a normal neighborhood of $x$.
Using normal coordinates and a local frame field of $E$,
introduce norms on $T_xX$ and uniformly on the fibers of $E$.
Write $c(t)=\exp_x w(t)$.
Define $\gamma_\eps(t)=\exp_x((t/\eps) w(\eps))$.
Observe that
\[ \sup\nolimits _{0\leq t\leq\eps} |w(t)-(t/\eps)w(\eps)|=\bigoh(\eps^2). \]
Standard Lipschitz stability estimates for systems of ordinary differential equations,
applied to \eqref{eq-ode-par-transp}, imply the estimate
$\|\PTE_{c|_{[0,\eps]}}-\PTE_{\gamma_\eps}\|=\bigoh(\eps^2)$ with respect to norms of $\Hom(E_x,E_{c(\eps)})$.
To complete the proof, observe that
$\PTE_{x\from c(\eps)}= \big(\PTE_{\gamma_\eps}\big)^{-1}$ stays uniformly bounded.
\end{proof}

Let $a\in\Cinfty(T^*X;\pi^* E)$.
The covariant derivative $\nabla^{\pi^* E}a$ is a section of the vector bundle $T^*(T^*X)\otimes \pi^*E$.
More useful are the vertical derivative $\vnabla a$ and the horizontal derivative $\hnabla a$
which are sections of the bundles $\pi^* (TX\otimes E)$ and $\pi^* (T^*X\otimes E)$, respectively.
As in the scalar case, the vertical derivative is the derivative in fiber direction:
$\vnabla a(x,\xi)\in T_xX\otimes E_x$ is the derivative of the map $T_x^*X\to E_x$, $\xi\mapsto a(x,\xi)$,
that is
$\vnabla a(x,\xi)=\partial_j\otimes \partial^j a(x,\xi)$
holds in canonical coordinates.
Viewing $\eta\in T_x^*X$ as a tangent vector to the fiber, $\eta\in T^*_{(x,\xi)}T^*X$,
we have $\vnabla a(x,\xi)\cdot\eta=\nabla^{\pi^*E}_\eta a(x,\xi)$.
We define the horizontal derivative with respect to a local frame $(e_j)$
with connection matrix $(\omega_k^j)$ as
\begin{equation}
\label{eq-horiz-a-vector}
\hnabla a= (\hnabla a^j + a^k \omega_k^j)e_j,
   \quad a(x,\xi)=a^j(x,\xi) e_j(x).
\end{equation}
For real-valued $\varphi\in\Cinfty(X)$ we have
\begin{equation}
\label{eq-covder-horiz-der}
\nabla^E(a\circ\intd\varphi)=(\hnabla a)\circ\intd\varphi + \nabla^2\varphi \cdot (\vnabla a\circ\intd\varphi).
\end{equation}
In particular, $\hnabla a(x,\xi)=\nabla^E (a\circ\intd\varphi)(x)$
if $\xi=\intd\varphi(x)$ and $\nabla^2 \varphi(x)=0$.
This shows that the horizontal derivative is well-defined.
Contractions of $T^*X\otimes TX$ are denoted by a centered dot, e.g.,
\[ \vnabla \cdot \hnabla a, \nabla^2\varphi\cdot\vnabla^2 a\in \Cinfty(T^*X;\pi^* E). \]
An example in canonical coordinates is
\[
\hnabla\cdot\vnabla a=\vnabla \cdot \hnabla a
   = \big(\partial^k\partial_k a^j + \Gamma_{km}^n (\delta^k_n \partial^m a^j +\xi_n \partial^k\partial^m a^j)
       + \partial^m a^k \omega^j_k(\partial_m)\big) e_j
\]
where $a=a^j e_j$; recall \eqref{eq-horiz-a-scalar}.
In the flat case, this simplifies to $\vnabla\cdot\hnabla a= \partial^k\partial_k a$,
which is a term familiar from the formula for the subprincipal symbol.
When dealing with products in the geometric symbol calculus we shall encounter the following situation:
$E$, $F$ and $G$ are vector bundles over $X$, $E$ and $F$ with connections.
The bundle $\Hom(E,F)$ is equipped with the induced connection.
Let $a$ and $b$ be sections of the bundles $\pi^* \Hom(E,F)$ and $\pi^* \Hom(F,G)$, respectively.
Besides $ba$ also
\(  \vnabla b \cdot \hnabla a \)
is a well-defined element of $\Cinfty(T^*X;\pi^* \Hom(E,G))$.

The Poisson bracket $\{b,a\}=H_b a$ and the Hamilton field $H_b$
are defined for functions $a,b\in\Cinfty(T^*X)$.
We have the following generalization of the Poisson bracket when $a$ is not necessarily scalar.
\begin{lemma}
Let $a\in\Cinfty(T^*X;\pi^*\End(E))$ and $b\in\Cinfty(T^*X)$ real-valued.
Denote by $H_b$ the Hamilton vector field of $b$.
Then
\begin{equation}
\label{eq-Poisson-vert-horiz}
\nabla^{\pi^*\End(E)}_{H_b} a = \vnabla b\Id\cdot\hnabla a - \hnabla b\Id\cdot\vnabla a.
\end{equation}
\end{lemma}
\begin{proof}
Write $a(x,\xi)=a^j(x,\xi) e_j(x)$.
Both sides of \eqref{eq-Poisson-vert-horiz} are equal to
\[
\big(\partial^k b (\partial_k a^j +\omega^{j}_i(\partial_k) a^i)
    -\partial_k b \partial^k a^j\big) e_j.
\]
To see that this is true for the right-hand side,
use \eqref{eq-horiz-a-vector}, \eqref{eq-horiz-a-scalar},
the symmetry of the Christoffel symbols, and $\nabla^{\pi^*\End(E)} \Id=0$.
On the left-hand side insert $H_b=(\partial^j b)\partial_j-(\partial_j b)\partial^j$.
\end{proof}

\section{Geometric symbol calculus}
\label{section-geom-symbol-calc}
Assuming additional geometric structure,
geometric pseudo-differential calculi invariantly define full symbols of operators;
see \cite{Bokobza69pdovardiff,Widom80completesymb} and the more recent
work of \cite{Sharaf04geosymb,Sharaf05geosymb}.
In this section we present,
following Sharafutdinov's approach and including complete proofs,
a geometric symbol calculus down to subprincipal (leading) symbol level.

Suppose $(X,\nabla)$ is a manifold with a symmetric connection,
and $E$ a complex vector bundle over $X$ endowed with a connection $\nabla^E$.
Let $u\in\Cinfty(X;E)$.
Using parallel transport, the covariant derivative is expressed as an ordinary derivative:
\[ \nabla^E u(x) = U'(0)\in\Hom(T_xX,E_x), \quad U(v) = \PTE_{x\from \exp_x v} u(\exp_x v).  \]
Suppose $\supp u$ is a compact subset of a normal neighborhood of $x$.
The Fourier inversion formula applied to the derivative $U'$ gives
\begin{equation}
\label{eq-nablaE-via-Fourier}
-i\nabla^E u(x)=
 (2\pi)^{-\dim X}\int_{T^*_x X}\int_{T_x X} e^{-i\xi v} \xi\otimes \PTE_{x\from y} u(y)\intd v \intd \xi,
  \quad y =\exp_x v.
\end{equation}
The Lebesgue measures $\intd v$ and $\intd\xi$ are normalized as follows:
$\intd v\intd \xi=|\sigma^n|/n!$
where $\sigma$ denotes the canonical symplectic form of $T^*_x X\times T_x X$.

Let $F$ be another complex vector bundle over $X$.
Let $p\in \Cinfty(T^*X;\Hom(\pi^*E,\pi^*F))$ belong to a standard symbol class $S^m$.
As in \cite[Lemma~8.1]{Sharaf05geosymb},
we define a geometric pseudo-differential operator $P$ by
\begin{equation}
\label{eq-def-geom-PsDO-pointwise}
Pu(x)= (2\pi)^{-\dim X}\int_{T^*_x X}\int_{T_x X} e^{-i\xi v} p(x,\xi) \PTE_{x\from y} u(y)\chi(x,y)\intd v \intd \xi \in F_x,
\quad y=\exp_x v,
\end{equation}
$u\in\Cinfty(X;E)$ and $x\in X$.
We require the cutoff function $\chi\in\Cinfty(X\times X)$ to satisfy the following:
\begin{compactenum}[(i)]
\item $\chi=1$ in a neighborhood of the diagonal.
\item\label{item-chi-proper}
The relation $\supp\chi$ is proper, i.e., for $K\subset X$ compact,
the intersection of $\supp\chi$ with $K\times X$ and with $X\times K$ is compact.
\item\label{item-chi-admissable}
For every point in $X$ there exist open neighborhoods $U_1$ and $U_2$, $U_2$ a normal neighborhood of $U_1$,
that is of every point of $U_1$,
such that $(U_1\times X)\cap\supp\chi\subset U_1\times U_2$.
\end{compactenum}
There exists an open neighborhood $U\subset X\times X$ of the diagonal
such that $\supp\chi\subset U$ implies condition \eqref{item-chi-admissable}.
For example, if $(O_j)$ is a locally finite covering of $X$ by normal convex open sets $O_j$,
then $U=\cup_j O_j\times O_j$ has this property.
Indeed, given $z\in X$, set $U_1=\cap_{j\in J}O_j$ and $U_2=\cup_{j\in J}O_j$,
where $J$ is the finite set of indices $j$ with $z\in O_j$.

The domain of integration in \eqref{eq-def-geom-PsDO-pointwise} depends on $x$.
It is not immediately clear that \eqref{eq-def-geom-PsDO-pointwise} defines
a continuous section $Pu$, let alone a pseudo-differential operator $P$.
\begin{lemma}
\label{lemma-geom-PsDO-local}
Let $U_1\subset X$ open, $U_2$ a normal neighborhood of $U_1$,
such that $(U_1\times X)\cap\supp\chi\subset U_1\times U_2$.
Let $x,z\in U_1$, $x=\exp_z s$.
Define $w\mapsto v$ by $\exp_x v=\exp_z w\in U_2$.
For $\zeta\in T_z^*X$, set $\xi =\mathbin{^t(\partial v/\partial w)^{-1}}\zeta\in T_xX$ and $\varphi =-\xi v$.
Then \eqref{eq-def-geom-PsDO-pointwise} becomes
\begin{equation}
\label{eq-def-geom-PsDO-local}
Pu(x)= (2\pi)^{-\dim X}\int_{T^*_z X}\int_{T_z X} e^{i\varphi} p(x,\xi) \PTE_{x\from y} u(y)\chi(x,y)\intd w \intd \zeta,
  \quad y=\exp_z w.
\end{equation}
The phase function satisfies $\varphi(s,w,\zeta)=\zeta\psi(s,w)(s-w)$
with $\psi$ a $\Cinfty$ map into $\End(T_zX)$, and $\psi(s,s)=\Id$.
Moreover, $\varphi'_{\zeta}=0$ iff $w=s$.
\end{lemma}
\begin{proof}
The map $(w,\zeta)\mapsto (v,\xi)$ is symplectic, hence volume preserving.
Changing variables from $(v,\xi)$ to $(w,\zeta)$ shows that the
integral \eqref{eq-def-geom-PsDO-pointwise} equals the integral \eqref{eq-def-geom-PsDO-local}.
In case the double integral is not absolutely convergent,
we employ the standard procedure of
partial integration in regions where the phase function non-stationary.

Note that $\varphi=-\zeta(\partial v/\partial w)^{-1} v$,
and $\varphi'_{\zeta}\neq 0$ iff $v\neq 0$.
The map $(s,w)\mapsto (\partial v/\partial w)^{-1}v$ is a smooth map into $T_zX$ which is zero when $w=s$.
Thus there exists a $\Cinfty$ map $\psi$ into $\End(T_zX)$ such that
$(\partial v/\partial w)^{-1} v= \psi(s,w)(w-s)$ holds.
Taking the derivative with respect to $w$ and evaluating at $w=s$, we see that $\psi(s,s)$ is the identity.
\end{proof}
In $U_1\times U_1$, the phase function $\varphi$ parametrizes the conormal bundle of the diagonal, which is given by $s=w$.
By Lemma~\eqref{lemma-geom-PsDO-local} and standard theory, we know
that \eqref{eq-def-geom-PsDO-pointwise} defines a pseudo-differential operator $P$.
Furthermore, it follows from \eqref{item-chi-proper} that $P$ is properly supported.
Replacing $\chi$ by a cutoff function, which also satisfies the assumptions,
modifies $P$ by a smoothing operator.
Conversely, if $P$ is a pseudo-differential operator, then any given point $z\in X$ has
an open neighborhood $U_1$, where $P$ can be represented, modulo a smoothing operator,
as \eqref{eq-def-geom-PsDO-local}.
Here, when starting from a matrix representing $P$ with respect to some local
frame of $E$, a parallel transport map is absorbed into the symbol $p$.
Reversing the transformation from \eqref{eq-def-geom-PsDO-pointwise} to \eqref{eq-def-geom-PsDO-local},
we see that every pseudo-differential operator is, modulo a smoothing operator, geometric.

Denote by $S^m(T^*X;\Hom(E,F))\subset\Cinfty(T^*X;\Hom(\pi^*E,\pi^*F))$
the standard space of symbols of order $\leq m$.
To ease writing, we omit $\pi^*$ from symbol space notation,
and often we abbreviate the symbol space as $S^m$.
Formula \eqref{eq-def-geom-PsDO-pointwise} defines the geometric quantization $p\mapsto P=\op(p)$. 
The sum of $\op(S^m(T^*X;\Hom(E,F)))$ and the space $\Psi^{-\infty}(X;E,F)$
of properly supported smoothing operators
is the space $\Psi^m(X;E,F)$ of pseudo-differential operators of order $\leq m$.
A symbol $p\in S^m$ is polyhomogeneous of degree $m$ iff there exist
an asymptotic symbol expansion $p\sim\sum_{j\geq 0} p_j$ where the $p_j$'s are
$\Cinfty$ sections over $T^*X\nozerosection$,
and $p_j(x,\xi)$ is homogeneous of degree $m-j$ in the fiber variable $\xi$.
Thus, using semiclassical notation with a positive asymptotic parameter $\hh$,
\begin{equation}
\label{eq-polyhom-expansion}
\hh^m p(x,\xi/\hh)\sim \sum\nolimits_{j\geq 0} \hh^j p_j(x,\xi) \quad\text{as $\hh\to 0+$.}
\end{equation}
$\Sphg^m\subset S^m$ and $\Psiphg^m\subset\Psi^m$ denote the subclasses of polyhomogeneous symbols and operators.
The left-hand side in \eqref{eq-polyhom-expansion} is the symbol 
of the semiclassical operator $\hh^mP$ associated with \eqref{eq-def-geom-PsDO-pointwise}:
\begin{equation}
\label{eq-P-semiclass}
\hh^mP u(x)
  = (2\pi \hh)^{-\dim X} \int_{T^*_{x} X} \int_{T_{x} X} e^{-i\xi v/\hh} \hh^mp(x,\xi/\hh)
          \PTE_{x\from y}u(y)\chi(x,y) \intd v \intd \xi,
\end{equation}
$y=\exp_x v$.
A linear operator $P$ is a (classical) pseudo-differential operator of order $\leq m$
iff a complete asymptotic expansion
$\hh^m e^{-i\varphi/\hh} P e^{i\varphi/\hh}a(x) \sim \sum\nolimits_{j\geq 0} \hh^j b_j(x)$
exists whenever the phase $\varphi$ is real-valued and $\intd\varphi\neq 0$ holds on $\supp a$.
We give a formula the two top order terms.
\begin{proposition}
\label{prop-fund-asymp-expansion}
Let $P=\op(p)\in\Psiphg^m(X;E,F)$, $p\sim\sum_{j\geq 0} p_j$.
Let $a\in\Cinfty(X;E)$, $\varphi\in\Cinfty(X)$ and realvalued, $\intd \varphi\neq 0$ on $\supp a$.
As $\hh\to 0+$,
\begin{equation}
\label{eq-fund-asymp-expansion}
\hh^m e^{-i\varphi/\hh} P e^{i\varphi/\hh}a
   = ((p_0+\hh p_1)\circ\intd\varphi) a  -i\hh D a + \bigoh(\hh^2)
\end{equation}
holds, where
\[ D a= (\vnabla p_0\circ\intd\varphi) \cdot\nabla^E a + \nabla^2\varphi\cdot(\vnabla^2 p_0\circ\intd\varphi) a/2. \]
If $E=F$ and the principal symbol is scalar, $p_0=qI$,
then $V=\vnabla q\circ\intd\varphi$ is a vector field, and
\[ D a= \nabla^E_V a + \nabla^2\varphi\cdot(\vnabla^2 p_0\circ\intd\varphi) a/2.  \]
Furthermore, if $P\in\Psi^{-\infty}$, then $p_j=0$ for all $j$.
\end{proposition}

By \cite[Theorems~7.7.5-6]{Hormander90anaOne}, the stationary phase formula
\begin{equation}
\label{eq-stationary-phase}
\int e^{i \varphi/\hh} a\intd x = {\det (H/2\pi i\hh)}^{-1/2}
 \sum\nolimits _{j<3N} (i\hh/2)^j \langle H^{-1} \partial,\partial\rangle^j
           \big(e^{i\rho/\hh}a\big)(0)/j! +\bigoh(\hh^{N})
\end{equation}
holds as $\hh\to 0+$, uniformly in parameters.
Here $\varphi$ is $\Cinfty$ and real-valued, $\varphi(0)=0$,
the origin $x=0$ is the only critical point of $\varphi$ in the support of $a\in\Ccinfty(\R^n)$,
and $H=\varphi''(0)$ is non-singular.
Furthermore, $\rho(x)=\varphi(x)-\langle Hx,x\rangle/2=\bigoh(|x|^3)$.
The following corollary of \eqref{eq-stationary-phase} will be useful.
Assume that the integration variable splits as $x=(y,\eta)\in\R^k\times\R^k$.
If $\varphi''_{\eta\eta}(0,0)=0$ and
if $\partial_y^\alpha\partial_\eta^\beta\rho(0,0)=0$ when $|\alpha|\leq 2$, then
\begin{equation}
\label{eq-stationary-phase-mod-hsquared}
\int e^{i \varphi(x)/\hh} a(x)\intd x 
 = |\det (\varphi''_{y\eta}/2\pi \hh)|^{-1} \big(a(0) + i\hh \langle H^{-1} \partial,\partial\rangle a(0)/2 +\bigoh(\hh^{2})\big).
\end{equation}
This follows by direct computation from \eqref{eq-stationary-phase} using $N = 2+k$.

\begin{proof}[Proof of Proposition~\ref{prop-fund-asymp-expansion}]
A uniform asymptotic expansion of $e^{-i\varphi/\hh} \hh^m P e^{i\varphi/\hh}a$ is known to exist.
So it suffices to evaluate the expansion at a given point $x$ which we hold fixed.
By \eqref{eq-P-semiclass},
\begin{equation}
\label{eq-A-with-osc-test}
(e^{-i\varphi/\hh} \hh^m P e^{i\varphi/\hh}a)(x)
  = (2\pi \hh)^{-\dim X} \int_{T^*_{x} X} \int_{T_{x} X} e^{i\Phi/\hh} \hh^m p(x,\xi/\hh)
          \PTE_{x\from y}a(y)\chi(x,y) \intd v \intd \xi.
\end{equation}
Here $\Phi=\Phi(v,\xi)=-\xi v+\tilde\varphi(v)-\varphi(x)$ with $\tilde\varphi(v)=\varphi(y)$ and $y=\exp_x v$.

The stationary point $(v,\xi)$ of $\Phi$ satisfies $v=0$
and $\xi=\tilde\varphi_v'(0)=\intd\varphi(x)$.
The Hessian is
\[ H=\begin{bmatrix} \tilde\varphi_{vv}''(0) & -I\\ -I& 0 \end{bmatrix}. \]
Apply \eqref{eq-stationary-phase-mod-hsquared} to \eqref{eq-A-with-osc-test}, and get
$e^{-i\varphi/\hh} \hh^mP e^{i\varphi/\hh}a = b_0+\hh b_1+\bigoh(\hh^2)$,
$b_0(x)= p_0(x,\xi)a(x)$, and
\[
b_1(x)= p_1(x,\xi)a(x)-i\vnabla p_0(x,\xi)\cdot \nabla^E a(x)
    -(i/2) \tilde\varphi_{vv}''(0)\cdot \vnabla^2 p_0(x,\xi)a(x).
\]
By the symmetry of $\nabla$, $\tilde\varphi_{vv}''(0) =\nabla^2\varphi(x)$.
Thus \eqref{eq-fund-asymp-expansion} holds.
In the scalar case, $p_0=q\Id$, the simpler formula for $Da$ is easily verified.

Fix $w\in E_x$ and $0\neq \eta\in T_x^*X$.
If $\varphi$ and $a$ satisfy $\varphi(\exp_x v)=\eta v$ and $\PTE_{x\from y}a(y)=w$ for $y$ near $x$,
the stationary phase formula \eqref{eq-stationary-phase} applied to \eqref{eq-A-with-osc-test} gives
\[
\hh^m (e^{-i\varphi/\hh} P e^{i\varphi/\hh}a)(x)
 \sim \sum\nolimits _{j\geq 0} (i\hh)^j \langle\partial_\xi,\partial_v\rangle^j \hh^m p(x,\xi/\hh) w|_{(v,\xi)=(0,\eta)} 
 = \hh^mp(x,\eta/\hh)w.
\]
Suppose $P\in\Psi^{-\infty}$.
It follows that $p(x,\eta/\hh)w =\bigoh(\hh^{\infty})$ holds, thus $p_j=0$ for all $j$.
\end{proof}

It follows from \eqref{eq-fund-asymp-expansion} that, modulo $\bigoh(\hh^2)$,
\[
(p_0+\hh p_1)(x,\xi) a(x) \equiv (e^{-i\varphi/\hh} \hh^m P e^{i\varphi/\hh}a)(x)
  \quad\text{if $\xi=\varphi'(x)$, $\nabla^2\varphi(x)=0=\nabla^E a(x)$.}
\]
We call $p_0+\hh p_1$ the \emph{leading symbol} of the operator $P$.
Even when $P$ is not semi-classical we stick to semi-classical notation for the leading symbol
in order to visibly distinguish $p_0$ from $p_1$.
The principal symbol $p_0$ does not depend on the connections $\nabla$ and $\nabla^E$, whereas $p_1$ does.
It follows from \eqref{eq-nablaE-via-Fourier} that the leading symbol of $-i\nabla^E$
is equal to its principal symbol $d:T^*X\to \Hom(TM,\Hom(TM,E))$ which is given by
$d(\xi):e\mapsto e\otimes\xi$ for $\xi\in T_x^*M$, $e\in E_x$.

Leading symbols refine the symbol calculus of pseudo-differential operators.

\begin{proposition}
\label{prop-geosymb-composition}
The product $QP$ of pseudo-differential operators $Q$ and $P$
with leading symbols $q+\hh q_1$ and $p+\hh p_1$
has leading symbol $qp+\hh (q_1p+qp_1 -i \vnabla q\cdot \hnabla p)$.
\end{proposition}
\begin{proof}
Denote the orders of $P$ and $Q$ by $m_P$ and $m_Q$.
Put $m=m_P+m_Q$.
Let $\varphi$, $a$ and $Da$ as in Proposition~\ref{prop-fund-asymp-expansion}.
It follows from \eqref{eq-fund-asymp-expansion} that the following
holds modulo $\bigoh(\hh^2)$:
\begin{align*}
\hh^m e^{-i\varphi/\hh} QP e^{i\varphi/\hh}a
   &\equiv \hh^{m_Q} e^{-i\varphi/\hh} Q e^{i\varphi/\hh}
           \big(((p+\hh p_1)\circ\intd\varphi) a -i\hh D a\big) \\
   &\equiv ((q+\hh q_1)\circ\intd\varphi) \big(((p+\hh p_1)\circ\intd\varphi) a -i\hh  D a\big) \\
   &\phantom{=}    -i\hh (\vnabla q\circ\intd\varphi) \cdot\nabla^F ((p\circ\intd\varphi) a)
       -i\hh  \nabla^2\varphi\cdot(\vnabla^2 q\circ\intd\varphi) (p\circ\intd\varphi)a/2.
\end{align*}
Fix $(x,\xi)\in T^*X\nozerosection$.
Choose $\varphi$ and $a$ such that $\xi=\intd\varphi(x)$, $\nabla^2\varphi(x)=0$, and $\nabla^E a(x)=0$ hold.
Then $Da(x)=0$, and, by Leibniz' rule and \eqref{eq-covder-horiz-der},
\[
\nabla^F ((p\circ\intd\varphi) a)(x)
  = \nabla^{\Hom(E,F)} (p\circ\intd\varphi)(x) a(x) +p(x,\xi) \nabla^E a(x)
  = \hnabla p(x,\xi)a(x).
\]
Summarizing, we have shown that
\[
e^{-i\varphi/\hh} \hh^m QP e^{i\varphi/\hh}a(x)
  \equiv (qp+\hh q_1p+\hh qp_1)(x,\xi)a(x) -i\hh (\vnabla q\cdot\hnabla p)(x,\xi)a(x)
\]
holds modulo $\bigoh(\hh^2)$, which completes the proof.
\end{proof}

\begin{remark}
\label{remark-mult-op}
Suppose $Q$ is a differential operator of order zero, that is,
the multiplication operator given by a bundle homomorphism $q$.
Then the leading symbol of $QP$ is $qp+\hh qp_1$.
\end{remark}

In the remainder of this section,
we assume that $X$ is an oriented (pseudo-)Riemannian manifold with metric $g$,
$\nabla$ the Levi-Civita covariant derivative,
and $\volg\in\Cinfty(X;\DB{X})$ the positive volume form.
Every $P\in\Psi^m(X;E,F)$ can be viewed as an operator acting on half-densities:
\[ \volg^{1/2}P\volg^{-1/2}\in\Psi^m(X;E\otimes\hDB{X},F\otimes\hDB{X}). \]
We identify symbol spaces via the identification
\[ \Hom(E,F)\equiv \Hom(E\otimes\hDB{X},F\otimes\hDB{X}) \]
which is defined by tensoring sections with $\volg^{\pm 1/2}$.

\begin{corollary}
\label{cor-leadsymb-inv-conjug}
The leading symbols of $P$ and $\volg^{1/2}P\volg^{-1/2}$ agree.
\end{corollary}
\begin{proof}
The vertical derivatives of the multiplication operators $\volg^{\pm 1/2}$ vanish.
So do the horizontal derivatives as $\volg$ is parallel, $\nabla\volg=0$.
Now apply Proposition~\ref{prop-geosymb-composition}.
\end{proof}

Next we consider formal adjoints of $P$ and of $Q=\volg^{1/2}P\volg^{-1/2}$.
Assume that $E$ and $F$ are Hermitian bundles with metric connections.
Denote the bundle metrics $(\cdot|\cdot)_E$ and $(\cdot|\cdot)_F$,
and the metric connections $\nabla^E$ and $\nabla^F$.
The adjoint of $p\in\Hom(E,F)$ is written $p^*\in\Hom(F,E)$.
For sections $u,v$ of $E$, resp.\ $E\otimes\hDB{X}$, we have
scalar products $\int_X(u|v)_E\volg$, resp.\ $\int_X(u|v)_E$.
Associated to the scalar products 
are the formal adjoints $P^*\in\Psi^m(X;F,E)$ and $Q^\dagger= \volg^{1/2}P^*\volg^{-1/2}$.

\begin{proposition}
\label{prop-geosymb-adjoint}
The leading symbol of the adjoint $P^*$ of a pseudo-differential operator $P$,
polyhomogeneous with leading symbol $p+\hh p_1$, equals $p^*+\hh p_1^* -i\hh  \vnabla \cdot \hnabla p^*$.
\end{proposition}
By Corollary~\ref{cor-leadsymb-inv-conjug}, the proposition holds also for $P^\dagger$.
\begin{proof}
Suppose $P\in\Psiphg^m(X;E,F)$.
Let $a\in\Ccinfty(X;E)$, $b\in\Cinfty(X;F)$, and
$\varphi\in\Cinfty(X)$ real-valued, $\intd\varphi\neq 0$.
Apply \eqref{eq-fund-asymp-expansion} to the right-hand side of
\begin{equation*}
\int_X \big(a\big|e^{-i\varphi/\hh}\hh^m P^* e^{i\varphi/\hh} b\big)_E\volg
  = \int_X \big(e^{-i\varphi/\hh}\hh^m P e^{i\varphi/\hh} a\big|b\big)_F\volg,
\end{equation*}
and get
\begin{align*}
\int_X \big(a\big|e^{-i\varphi/\hh}\hh^m P^* e^{i\varphi/\hh} b\big)_E\volg
   &= \int_X \big(a\big|((p^*+\hh p_1^*)\circ\intd\varphi+(i\hh/2)\nabla^2\varphi\cdot\vnabla^2 p^*\circ\intd\varphi)b\big)_E\volg \\
   &\phantom{==} -i\hh \int_X\big((\vnabla p\circ\intd\varphi)\cdot \nabla^E a\big| b\big)_F\volg
         +\bigoh(\hh^2)
\end{align*}
as $\hh\to 0+$.
Suppose $(V_j)$ is a frame of $TX$ over an open subset $O\subset X$, and $\supp b\subset O$.
Write $\vnabla p^*\circ\intd\varphi=\sum_j S^j\otimes V_j$ with $S^j\in\Cinfty(O;\Hom(F,E))$.
Then
\[
\big((\vnabla p\circ\intd\varphi)\cdot\nabla^E a\big|b\big)_F
  = \sum\nolimits _j \big(\nabla^E_{V_j}a\big| S^j b\big)_E
  = -\sum\nolimits _j \big(a\big|\nabla^E_{V_j} S^j b\big)_E
    +\sum\nolimits _j V_j \big(a\big| S^j b\big)_E.
\]
By the Leibniz' rule,
\begin{align*}
\nabla^E_{V_j} S^jb 
   &= (\nabla^{\Hom(F,E)}_{V_j} S^j) b + S^j \nabla^F_{V_j} b, \\
\nabla^{\Hom(F,E)\otimes TX}(S^j\otimes V_j)
   &= \nabla^{\Hom(F,E)} S^j \otimes V_j + S^j\otimes \nabla V_j.
\end{align*}
Using \eqref{eq-covder-horiz-der} and contraction, we derive
\[
\hnabla\cdot\vnabla p^* \circ \intd\varphi + \nabla^2\varphi\cdot \vnabla^2 p^*\circ \intd\varphi
  = \sum\nolimits_j\nabla^{\Hom(E,F)}_{V_j} S^j +\sum\nolimits_j (\divergence V_j)S^j.
\]
Combining formulas, we obtain
\begin{align*}
\big((\vnabla p\circ\intd\varphi)\cdot \nabla^E a\big| b\big)_F
   &= -\big(a\big|(\hnabla\cdot\vnabla p^* \circ \intd\varphi + \nabla^2\varphi\cdot \vnabla^2 p^*\circ \intd\varphi) b\big) \\
   &\phantom{==} + \sum\nolimits_j (V_j c^j +(\divergence V_j)c^j)
                 - \sum\nolimits_j (a|S^j \nabla^F_{V_j} b)_E,
\end{align*}
$c^j=(a|S^j b)_E$.
Integrate and use \eqref{eq-int-LieV}.
Since $a$ is arbitrary, we deduce the asymptotics
\begin{align*}
e^{-i\varphi/\hh}\hh^m P^* e^{i\varphi/\hh} b
  &= ((p^*+\hh p_1^*)\circ\intd\varphi ) b -i\hh (\hnabla\cdot\vnabla p^* \circ \intd\varphi) b \\
  &\phantom{==} -(i\hh/2)\nabla^2\varphi\cdot\vnabla^2 p^*\circ\intd\varphi)b 
   -i\hh S^j \nabla^F_{V_j} b +\bigoh(\hh^2).
\end{align*}
Fix $x\in X$ and $\xi=\intd\varphi(x)$, and suppose that $\nabla^2\varphi(x)=0$ and $\nabla^F b(x)=0$. 
Then
\[
(e^{-i\varphi/\hh}\hh^m P^* e^{i\varphi/\hh} b)(x)
   = (p^*+\hh p_1^*)(x,\xi)b(x) -i\hh (\hnabla\cdot\vnabla p^*)(x,\xi)b(x) + \bigoh(\hh^2),
\]
which implies the proposition.
\end{proof}

We apply the geometric calculus to the elasticity operator
$L=\nabla^*C\nabla=(-i\nabla)^*C(-i\nabla)$ of an elastic body $(M,g,C)$.
The leading symbol $d(\xi)\in \Hom(T_xM,T_x^{(1,1)}M)$ of $-i\nabla$
at $\xi\in T_x^*M$ is given by tensoring with $\xi$.
The horizontal derivatives of $d$ and of its adjoint symbol $d^*$ vanish.
The adjoint $(-i\nabla)^*$ is defined by taking the $L^2$ inner product given by $\intd V_M=\volg$.
Proposition~\ref{prop-geosymb-adjoint} implies that $d^*(\xi)$, which is contraction with $\xi$,
is the leading symbol of $(-i\nabla)^*$.
Apply Remark~\ref{remark-mult-op} to multiplication by the stiffness tensor $C$.
It follows from Proposition~\ref{prop-geosymb-composition} that the leading symbol of $L$ equals
\[
\ell+\hh\ell_1, \quad
\ell(\xi) = d(\xi)^* C(x) d(\xi), \quad
\ell_1(\xi) = -i \big((\vnabla d(\xi)^*)\cdot \hnabla C\big) d(\xi),
\]
where $\xi\in T_x^*M$, and $d(\xi):w\mapsto w\otimes\xi$, $w\in T_xM$.
Given local coordinates $x^j$, this reads
\[ \ell(\xi)=\big(C^{ijkm}\xi_j\xi_m\big), \quad \sqrt{-1}\ell_1(\xi)=\big(C^{ijkm}_{\phantom{ijkm};m}\xi_j\big). \]
The $\xi_j$'s are dual coordinates,
and the matrix representations correspond to the frame $\partial_j$.
Thus the lower order symbol $\sqrt{-1}\ell_1$ is a contraction of the tensor $\nabla C$.
Using the Leibniz rule and the symmetries of $C$, we find the formula
\begin{equation}
\label{eq-ell-leadsymb-simplified}
2i\ell_1 =\vnabla\cdot\hnabla \ell,
\end{equation}
$i=\sqrt{-1}$.
This means that the leading symbol of the elasticity operator is determined from its principal symbol.

\section{Subprincipal symbol}

In this section we define the subprincipal symbol of systems of real principal type.
Examples from elastodynamics are considered.

Let $P\in\Psiphg^m(X;E)$ with principal symbol $p$ and
characteristic $\Char P=\{\det p=0\}\nozerosection$.
Following \cite{Dencker82polarrpt}, $P$ is of real principal type iff there hold:
\begin{compactenum}[(i)]
\item
\label{rpt-Hamilton}
$\Char P$ is smooth hypersurface with a non-radial Hamilton field $H$
which is homogeneous of degree zero.
\item
Locally on $\Char P$, the dimension of $\kernel p$, the null-space of $p$, is constant.
\item
\label{rpt-tilde-p}
If $q$ is a defining function for $\Char P$,
then $\tilde p=qp^{-1}$ extends smoothly across $\Char P$.
\end{compactenum}
Fixing $H$ as in \eqref{rpt-Hamilton}, and we say that $P$ is endowed with the Hamiltion field $H$.
So the orientation of the bicharacteristic strips, that is, the integral curves of $H$ in $\Char P$, is determined.
Furthermore, up to a factor which equals $1$ on $\Char P$, Hamilton functions $q$
satisfying $H=H_q$ on $\Char P$ are uniquely determined.
So is $\tilde p|_{\Char P}$.
The range of $\tilde p$ equals $\kernel p$, and the range of $p$ equals $\kernel \tilde p$.
The definition of real principal type systems microlocalizes to open conic subsets of $T^*X\nozerosection$.

\begin{proposition}
Suppose $P$ is of real principal type and endowed with the Hamilton field $H=H_q$.
Denote by $p+\hh p_1$ the leading symbol of $P$.
Set
\begin{equation}
\label{eq-subprinc-symb}
p_s =\tilde p p_1 -i \vnabla\tilde p\cdot \hnabla p +(i/2) (\vnabla\cdot\hnabla q)\Id,
\end{equation}
where $q\Id=\tilde p p$.
The restriction $p_s|_{\kernel p}$ does not depend on the choice of $q$.
\end{proposition}
We call $p_s$, and its restriction to the kernel of $p$,
the \emph{subprincipal symbol} of $P$ with respect to the Hamilton field $H$.
Note that $p_s$ is homogeneous of degree zero because we do assume that $q$ is homogeneous of degree $1$.
For scalar operators of order $1$, our subprincipal symbol agrees with the usual one.
\begin{proof}
Suppose $H=H_q=H_{q'}$.
Thus $q'=cq$ with $c=1$ at $q=0$.
Moreover, $q=\tilde p p$ and $q'=\tilde p' p$, $\tilde p'=c\tilde p$.
We identify scalar functions $f$ with the corresponding sections $f\Id$.
Since $\Id=\Id_E$ is parallel, \eqref{eq-Poisson-vert-horiz} implies
$\vnabla q\cdot\hnabla c-\vnabla c\cdot\hnabla q=(Hc)\Id=0$.
Therefore, at $q=0$,
\[ \vnabla\cdot\hnabla q'
  = \vnabla\cdot\hnabla q+\vnabla c\cdot\hnabla q +\vnabla q\cdot\hnabla c
  = \vnabla\cdot\hnabla q+2\vnabla c\cdot\hnabla q.
\]
Using $\hnabla q=(\hnabla\tilde p)p+\tilde p(\hnabla p)$, we have
\(
\vnabla\tilde p'\cdot \hnabla p
  = c\vnabla\tilde p\cdot \hnabla p + \vnabla c\cdot\big(\hnabla q -(\hnabla\tilde p) p\big).
\)
It follows that
\[
\vnabla\tilde p'\cdot \hnabla p -\vnabla\cdot\hnabla q'/2
 = \vnabla\tilde p\cdot \hnabla p -\vnabla\cdot\hnabla q/2
  -(\vnabla c\cdot \hnabla\tilde p) p
\]
holds at $q=0$.
Therefore, $p_s|_{\kernel p}$ does not depend on $c$.
\end{proof}

Suppose $\Lambda$ is a Lagrangian submanifold of $T^*X$ which is contained in $\Char P$.
Suppose $E$ equipped with a connection $\nabla^E$.
Denote by $(\hat E,\nabla^{\hat E})$ the pullback of $(E,\nabla^E)$ by the canonical projection $\Lambda\to X$,
and by $\hDB{\Lambda}$ the bundle of half-densities over $\Lambda$.
Setting
\[ D_H (e\otimes\nu) = (\nabla^{\hat E}_H e)\otimes\nu + e\otimes(\Lie_H \nu), \]
the Hamilton field $H$ and the subprincipal symbol $p_s$ determine a first order differential
operator $D_H$ which operates on sections of the bundle $\hat E\otimes\hDB{\Lambda}\to \Lambda$.
Here $e$ is a section of $\hat E$,
and $\Lie_H \nu$ denotes the Lie derivative of a nowhere vanishing the half-density $\nu$.
Note that $D_H$ is well-defined.

\begin{lemma}
\label{lemma-homog-transp-eq}
Suppose $a\in\Cinfty(\Lambda;\hat E\otimes\hDB{X})$ satisfies $i^{-1}D_H a +p_sa\in\kernel p$.
Then $pa$ vanishes along a bicharacteristic strip of $H$ if $pa$ vanishes at some point on it.
\end{lemma}
\begin{proof}
$\vnabla q=(\vnabla p)\tilde p +p\vnabla\tilde p$
and $\hnabla q=(\hnabla\tilde p)p+\tilde p\hnabla p$
imply the first equality of
\[
p\vnabla\tilde p\cdot \hnabla p -\vnabla p\cdot(\hnabla \tilde p)p
   = \vnabla q\cdot \hnabla p- \vnabla p\cdot \hnabla q
   = \nabla_H^{\pi^*\End(E)}p.
\]
The second follows from \eqref{eq-Poisson-vert-horiz}.
Observe that $D_H pa= (\nabla_H^{\End(\hat E)} p)a+ pD_Ha$, and that
$\nabla_H^{\End(\hat E)} p$ equals the pullback of $\nabla_H^{\pi^*\End(E)} p$ to $\Lambda$.
Using the assumption and \eqref{eq-subprinc-symb}, we get
\begin{align*}
D_H pa
  &=\big(p\vnabla\tilde p\cdot \hnabla p -\vnabla p\cdot(\hnabla \tilde p)p-ipp_s\big)a \\
  &= -\vnabla p\cdot(\hnabla \tilde p)pa +(\vnabla\cdot\hnabla q)pa/2.
\end{align*}
This means that the restriction of $pa$ to a bicharacteristic satisfies
a homogeneous system of ordinary differential equations.
The assertion follows from uniqueness of solutions to initial value problems.
\end{proof}

Consider the operator $P=L-\rho D_t^2$ of elastodynamics for a three-dimensional elastic body.
The principal symbol is $p(\xi,\tau)=\ell(\xi)-\rho(x)\tau^2$ at $(\xi,\tau)\in T_{(x,t)}^*(M\times \R)$.
Suppose the stiffness tensor is isotropic with Lam\'e parameters satisfying $\lambda+\mu>0$ and $\mu\geq 0$.
The principal symbol $\ell(\xi)$ of the isotropic elasticity operator $L$ is stated
in formula \eqref{eladyn-iso-acoustic-tensor} of Appendix~\ref{sect-elastic-symmetries}.
Therefore, the principal symbol $p$ of the operator of elastodynamics for isotropic media is given by
\begin{equation}
\label{eq-p-iso-eladyn}
p(\xi,\tau)/\rho = (c_p^2\xi^2-\tau^2)\pi_p(\xi) + (c_s^2\xi^2-\tau^2)\pi_s(\xi).
\end{equation}
Here $c_s=\sqrt{\mu/\rho}<\sqrt{(\lambda+2\mu)/\rho}=c_p$ are the speeds of shear and pressure waves.
If $(\xi,\tau)\in\Char P$, then one of $\tau=\pm c_{p/s}|\xi|$ holds.
Off $\Char P$ there holds
\begin{equation*}
\rho p(\xi,\tau)^{-1} = (c_p^2\xi^2-\tau^2)^{-1}\pi_p(\xi) + (c_s^2\xi^2-\tau^2)^{-1}\pi_s(\xi).
\end{equation*}
So we see that $P$ is of real principal type.
For a fluid, $\mu=0$, $P$ is of real principal type where $\tau\neq 0$, despite $c_s$ being zero.

Now suppose the elastic medium is non-isotropic and that the eigenvalues of
the acoustic tensor $\ell(\xi)$, $\xi\neq 0$, are pairwise distinct.
Then $P$ is of real principal type.
In fact, setting $\tilde p=(\partial_\tau q) (\partial_\tau p)^{-1}=-(\partial_\tau q/2)\Id$ on $\kernel p$,
the extension required in \eqref{rpt-tilde-p} is seen to exists, \cite[Proposition 3.2]{Dencker82polarrpt}.
By Proposition~\ref{prop-transv-iso-perturb}, this observation applies to transversely isotropic media
which are generic small perturbations from isotropy.
The axis of rotational symmetry, $J(x)$, should not be parallel to the propagation direction $\xi$.

Suppose the elastodynamics operator $P=L-\rho D_t^2$ is of real principal type.
A natural endowment of $P$ with a Hamilton field is given by the condition $Ht=1$,
so that bicharacteristic strips are parametrized by time.
The leading symbol $p+\hh p_1$ is given by $p=\ell-\rho \tau^2$ and $p_1=\ell_1$.
Here we equipped $X=M\times\R$ with the connection induced from $M$ by projection off the time axis.
The projection is also used to pull back the (complexified) tangent bundle of $M$ to get the bundle $E\to X$.
The formula \eqref{eq-subprinc-symb} for the subprincipal symbol can be simplified to a formal Poisson bracket:
\begin{equation}
\label{eq-subprinc-eladyn}
2ip_s = \vnabla\tilde p\cdot \hnabla p - \hnabla\tilde p \cdot\vnabla p \quad\text{on $\kernel p$.}
\end{equation}
In fact, using \eqref{eq-ell-leadsymb-simplified}, \eqref{eq-subprinc-symb}, and $\partial_t p=0$, we compute
\begin{align*}
2ip_s &= \tilde p \vnabla\cdot\hnabla p + 2\vnabla\tilde p\cdot \hnabla p -\vnabla\cdot\hnabla \tilde p p \\
      &= \vnabla\tilde p\cdot \hnabla p -\hnabla \tilde p\cdot\vnabla p - (\vnabla\cdot\hnabla \tilde p)p,
\end{align*}
and this implies \eqref{eq-subprinc-eladyn}.

\section{Lagrangian solutions}
\label{section-Lagr-solns}

With a real-valued non-degenerate phase function $\varphi$ and a symbol $a$
there is associated the oscillatory integral
$\int e^{i\varphi(x,\theta)}a(x,\theta)\intd\theta$.
We assume that $a$ is a section of the bundle $E$ pulled back by the
map $(x,\theta)\mapsto x$, i.e., $a(x,\theta)\in E_x$ holds.
The duality pairing with test functions in $\Ccinfty(X;E^*\otimes\DB{X})$,
$E^*$ the dual bundle of $E$, is defined.
The oscillatory integral is a distribution section of $E$.
The wavefront set is contained in the Lagrangian manifold parametrized by $\varphi$.
Lagrangian distributions are microlocally sums of oscillatory integrals.
Let $\Lambda\subset T^*X\nozerosection$ be a closed conic Lagrangian submanifold.
Denote by $I^\mu(X,\Lambda;E)$ the space of Lagrangian distribution sections of $E$
of order $\leq \mu$ which are associated with $\Lambda$.
Denote the Maslov bundle by $M_\Lambda$
and the pullback of $E$ under the projection $\Lambda\to X$ by $\hat E$.
There is a principal symbol isomorphism
\[
I^\mu(X,\Lambda;\hDB{X}\otimes E)/I^{\mu-1} \equiv
S^{\mu+n/4}(\Lambda;M_\Lambda\otimes \hDB{\Lambda}\otimes\hat E)/S^{\mu-1+n/4},
\]
\cite[Theorem 25.1.9]{Hormander85anaFour}.
We assume that principal symbols are homogeneous.
The principal symbol $a$ of $A\in I^\mu(X,\Lambda;\hDB{X}\otimes E)$ at $(y,\eta)\in\Lambda$ is given by
\begin{equation}
\label{eq-symbol-Lagrdistr}
\langle\chi(y), a(y,\eta;\lambda)\rangle \omega^{-1/2}
   =\lim_{\hh\to 0+} (2\pi)^{-n/4} \hh^{\mu-n/4} \int_X e^{-i\psi/\hh} \langle\chi, A\rangle,
\end{equation}
$n=\dim X$;
\cite[\S 4.1]{Duistermaat73fio}.
Here the real-valued function $\psi\in\Cinfty(X)$ satisfies $\psi(y)=0$ and $\intd \psi(y)=\eta$.
The tangent plane at $(y,\eta)$ of the graph of $\intd\psi$ is denoted $\lambda$,
and it is assumed that $\lambda$ is transversal to $T_{(y,\eta)}\Lambda$.
The section $\chi\in\Ccinfty(X;E^*\otimes\hDB{X})$ vanishes outside a small neighborhood of $y$.
Finally, $\omega$ is the canonical volume form of $T^*X$,
and angular brackets are duality brackets.

A convenient representation of the inverse of the principal symbol map is found microlocally using
canonical coordinates $x,\xi$ such that the planes $\xi=$constant are transversal to $\Lambda$.
Then $\Lambda$ is parametrized by a phase function $\langle\xi,x\rangle-\Theta(\xi)$,
$\Theta\in\Cinfty$ real-valued and homogeneous of degree one.
Thus $\Lambda$ is given by $x=\Theta'(\xi)$, \cite[Theorem~21.2.16]{Hormander85anaThree}.
If $a(\xi)|\intd\xi|^{1/2}$ is the principal symbol of $A\in I^\mu$, then
\begin{equation*}
A(x) \equiv (2\pi)^{-3n/4}
      \int_{\R^n}e^{i(\langle\xi,x\rangle-\Theta(\xi))} \PTE_{x\from \Theta'(\xi)}a(\xi)|\intd x|^{1/2}\intd\xi
\quad \mod I^{\mu-1}.
\end{equation*}
To see this, substitute $\hh\xi$ for $\xi$, choose $\psi(x)=\langle\eta,x\rangle$,
and evaluate \eqref{eq-symbol-Lagrdistr} by stationary phase.

We generalize the Duistermaat--Hörmander formula \cite[Theorem~25.2.4]{Hormander85anaFour}
about products with vanishing principal symbol to non-scalar real principal type systems.

\begin{theorem}
\label{theorem-amp-transport-eqn}
Suppose $P\in\Psiphg^m(X;E)$ is of real principal type with leading symbol $p+\hh p_1$
and endowed with the Hamilton field $H$.
Let $\Lambda$ be a closed conic Lagrangian submanifold of $\Char P$.
Let $A\in I^\mu(X,\Lambda;\hDB{X}\otimes E)$ with the principal symbol $a$.
Suppose $PA\in I^{\mu+m-1}(X,\Lambda;\hDB{X}\otimes E)$.
Then $p a=0$, and the principal symbol $b$ of $PA$ satisfies
\begin{equation}
\label{eq-amp-transport}
i^{-1} D_H a + p_s a = \tilde p b.
\end{equation}
\end{theorem}
It follows from \eqref{eq-amp-transport} and $\kernel\tilde p=\image p$
that $b$ is uniquely determined modulo the range of $p$.
Maslov factors do not appear in \eqref{eq-amp-transport} because they are locally constant.
A non-geometric coordinate version of \eqref{eq-amp-transport}
is given in \cite[Theorem 3.1]{SHRoehrig04lagrsol}.
\begin{proof}
We work microlocally near a given point in $\Lambda$.
Choose $\tilde P\in\Psiphg^{-m+1}(X;E)$ with principal symbol $\tilde p$.
The product $Q=\tilde P P$ is of real principal type.
By Proposition \ref{prop-geosymb-composition}, the leading symbol $q_0+\hh q_1$ satisfies
$q_0=q\Id$ and, on $\kernel p$, $q_1=\tilde p p_1 -i \vnabla \tilde p\cdot \hnabla p$.
Therefore, equipping $Q$ with the Hamilton field $H$,
the subprincipal symbol of $Q$ is equal to $p_s$.

It suffices to prove the theorem with $P$ replaced by $Q$.
So we assume $m=1$, $p=q\Id$, $\tilde p=\Id$, and $p_s=p_1+(i/2)(\vnabla\cdot\hnabla q)\Id$.
Fixing canonical coordinates as above, $\Lambda=\{x=\Theta'(\xi)\}$ holds.
Recall the abbreviations $\partial_j=\partial/\partial x^j$ and $\partial^j=\partial/\partial \xi_j$.
Using Taylor's formula, we write $q(x,\xi)=q_j(x,\xi)(x^j-\partial^j\Theta(\xi))$.
(Summation over indices in opposite position is implied.)
On $\Lambda$, the Hamilton field reads $H= - q_j\partial^j$.

Put $\tilde a(x,\xi)=\PTE_{x\from \Theta'(\xi)}a(\xi)\in E_x$.
Suppose
\begin{equation}
\label{eq-A-Lagr-in-adapted-coords}
A(x) = (2\pi)^{-3n/4} \hh^{-\mu-3n/4}
      \int_{\R^n}e^{i(\langle\xi,x\rangle-\Theta(\xi))/\hh} \tilde a(x,\xi)|\intd x|^{1/2}\intd\xi.
\end{equation}
The principal symbol of $A$ is $a(\xi)|\intd\xi|^{1/2}$.
Apply Proposition~\ref{prop-fund-asymp-expansion} to the oscillatory integral \eqref{eq-A-Lagr-in-adapted-coords}, and get
\begin{align*}
PA(x) &= (2\pi)^{-3n/4}\hh^{-\mu-m-3n/4}
      \int e^{i(\langle\xi,x\rangle-\Theta(\xi))/\hh} \tilde b(x,\xi,\hh)|\intd x|^{1/2}\intd\xi, \\
\tilde b(x,\xi,\hh) &= q(x,\xi)\tilde a(x,\xi)+\hh p_1(x,\xi)\tilde a(x,\xi) 
     -i\hh\nabla_V^E \tilde a(x,\xi) \\
  &\phantom{==} -i\hh \big((\nabla_V|\intd x|^{1/2})|\intd x|^{-1/2}
                +\nabla^2\langle\xi,x\rangle\cdot \vnabla^2 q(x,\xi)/2\big)\tilde a(x,\xi)
  +\bigoh(\hh^2)
\end{align*}
with $V= (\partial^j q) \partial_j$.
A partial integration proves
\[
\int e^{i(\langle\xi,x\rangle-\Theta(\xi))/\hh} q(x,\xi)\tilde a(x,\xi)\intd \xi
   = i\hh\int e^{i(\langle\xi,x\rangle-\Theta(\xi))/\hh} \partial^j q_j(x,\xi) \tilde a(x,\xi)\intd \xi.
\]
Hence
$PA(x) \equiv (2\pi)^{-3n/4} \int e^{i(\langle\xi,x\rangle-\Theta(\xi))} b(x,\xi)|\intd x|^{1/2}\intd\xi$
modulo $I^{\mu+m-2}$, where
\[
b =  i\partial^j (q_j\tilde a)+p_1\tilde a -i\nabla_V^E \tilde a 
     +i(\Gamma^k_{jk}\partial^j q + \Gamma^j_{k\ell}\xi_j \partial^k\partial^\ell q)\tilde a/2.
\]
Here the $\Gamma$'s are the Christoffel symbols of the symmetric connection $\nabla$,
and we used the formulas $\nabla|\intd x|^{1/2}=-\Gamma^k_{jk}\intd x^j|\intd x|^{1/2}/2$
and $\nabla^2\langle\xi,x\rangle= \nabla \xi_j\intd x^j= - \Gamma^j_{k\ell}\xi_j \intd x^k\intd x^\ell$.
The principal symbol of $PA$ is $b(\Theta'(x),\xi)|\intd\xi|^{1/2}$.
Observe that $\nabla^E \tilde a(x,\xi)=0$ holds if $x=\Theta'(\xi)$.
We claim that, at $\Lambda$, the following holds:
\begin{equation}
\label{eq-claim-DH-a-mu}
D_H(a\nu) = -\partial^j (q_j\tilde a)\nu +(\partial^j\partial_j q) a\nu/2,
   \quad \nu=|\intd\xi|^{1/2}.
\end{equation}
Assuming \eqref{eq-claim-DH-a-mu} we obtain
\[
b\nu = i^{-1} D_H (a\nu) + p_1 a\nu 
     +i\big(\partial^j\partial_jq+\Gamma^k_{jk}\partial^j q + \Gamma^j_{k\ell}\xi_j \partial^k\partial^\ell q\big)a\nu/2
  \quad\text{at $\Lambda$.}
\]
Recall \eqref{eq-horiz-a-scalar}.
The sum in brackets expresses $\vnabla\cdot\hnabla q$ in coordinates, implying \eqref{eq-amp-transport}.

It remains to prove \eqref{eq-claim-DH-a-mu}.
Let $(y(t),\eta(t))$ denote the bicharacteristic curve in $\Lambda$
which passes at $t=0$ through $(x,\xi)$, $x=\Theta'(\xi)$.
This means that
\[ \frac{\intd}{\intd t} \eta_j(t) =-q_j(y(t),\eta(t)),\quad \eta(0)=\xi,\quad y(t)=\Theta'(\eta(t)).  \]
Set $\beta_t=y|_{[0,t]}$.
Denote by $\lambda_t$ the loop which consists of $\beta_t$ followed by the geodesic from $y(t)$ to $y(0)=x$.
There holds
\[
-q_j\partial^j \tilde a|_{x=\Theta'(\xi)}
   = \frac{\intd}{\intd t} \tilde a(x,\eta(t))\big|_{t=0}
   = \frac{\intd}{\intd t} \PTE_{\lambda_t} \PTE_{\beta_t^{-1}}  a(\eta(t))\big|_{t=0}
   = \nabla_H^{\pi^* E} a(\xi).
\]
The last equality follows from Lemma~\ref{lemma-trivial-holonomy} and the Leibniz rule.
Moreover, since $H$ is tangent to $\Lambda$, we may replace $\pi^* E$ by $\hat E$.
Note that $\partial^j\partial_j q = \partial^j q_j - (\partial_jq_k)\partial^j\partial^k\Theta$ holds on $\Lambda$.
The Lie derivative of the half-density $\nu$ is obtained as follows:
\[
(\Lie_H\nu)/\nu
  = -\partial^j q_j(\Theta'(\xi),\xi) /2
  = -\big((\partial_k q_j)\partial^k\partial^j \Theta +(\partial^j q_j)\big) /2
  = \partial^j\partial_j q/2 -\partial^j q_j.
\]
Compare \cite[(25.2.11)]{Hormander85anaFour}.
Combining the results proves \eqref{eq-claim-DH-a-mu}.
\end{proof}

\begin{remark}
Given $B\in I^{\mu+m-1}=I^{\mu+m-1}(X,\Lambda;\hDB{X}\otimes E)$, we solve $PA\equiv B$
modulo $I^{\mu+m-2}$ as follows.
First we solve transport equations \eqref{eq-amp-transport} to find $A\in I^\mu$ such that $pa=0$,
and such that the principal symbol of $B-PA\in I^{\mu+m-1}$ lies in $\kernel \tilde p$.
Using $\kernel \tilde p=\image p$, we then replace $A$ by $A+A'$,
where $A'\in I^{\mu-1}$ solves $PA'\equiv B-PA$ modulo $I^{\mu+m-2}$.
\end{remark}

We use the Lagrangian intersection calculus of \cite{MelroseUhlmann79intersection}
to solve real principal type equations with Lagrangian sources.
Suppose $(\Lambda_0,\Lambda)$ is a pair of conic Lagrangian submanifolds of $T^*X\nozerosection$
intersecting cleanly in the boundary of $\Lambda$, $\partial\Lambda=\Lambda\cap\Lambda_0$.
Let $A$ be an element of 
\[ I^{k}(\Lambda_0\cup\Lambda) = I^{k}(X,\Lambda_0\cup\Lambda;\hDB{X}\otimes E), \]
the space of Lagrangian distributions of order $\leq k$ associated with the intersecting pair.
To simplify formulas, we denote intersecting pairs by their union,
and we often drop $X$ and the bundle $\hDB{X}\otimes E$ from the notation.
There holds $\WF A\subset\Lambda_0\cup\Lambda$.
In $\Lambda\setminus\Lambda_0$, resp.\ in $\Lambda_0\setminus\Lambda$,
$A$ belongs to $I^k(\Lambda)$ with principal symbol $a$,
resp.\ to $I^{k-1/2}(\Lambda_0)$ with principal symbol $a_0$.
The principal symbol $(a_0,a)$ of $A$ satisfies a compatibility condition at $\partial\Lambda$ which we recall.
Suppose $\gamma\in\partial\Lambda$.
Choose smooth functions $s_1,\ldots,s_{n-1}$ and $q,t$ such that the following hold.
On $\partial\Lambda$, $\sigma=\intd s_1\wedge \ldots\wedge \intd s_{n-1}\neq 0$ at $\gamma$.
The function $q$ vanishes on $\Lambda$ and its differential on $\Lambda_0$ at $\gamma$ is nonzero.
The function $t$ vanishes on $\Lambda_0$, is positive on $\Lambda\setminus\Lambda_0$, and $H_qt>0$ at $\gamma$.
The Hamilton fields of $q$ and $t$ are $H_q$ and $H_t$.
Choose a Lagrangian plane $\mu\subset T_\gamma(T^*X)$ which is transversal to the fiber
and to the connecting path between $T_\gamma \Lambda_0$ and $T_\gamma \Lambda$
generated from $T_\gamma \partial\Lambda$ and the straight line between $H_t$ and $H_q$.
The symbols $a$ and $qa_0$ are smooth at $\partial\Lambda$, and
\begin{equation}
\label{eq-iLagr-compat}
a|\sigma\wedge \intd t|^{-1/2}
  = (2\pi)^{1/4}e^{\pi i/4} q a_0|\sigma\wedge \intd q|^{-1/2} (H_qt)^{1/2}
\end{equation}
holds at $\gamma$ and $\mu$.
Conversely, if $(a_0,a)\in S^{k+(n-2)/4}(\Lambda_0)\times S^{k+n/4}(\Lambda)$
satisfy the compatibility condition \eqref{eq-iLagr-compat},
then there exists $A\in I^k(\Lambda_0\cup\Lambda)$, unique modulo $I^{k-1}(\Lambda_0\cup\Lambda)$,
which has $(a_0,a)$ as its principal symbol.
On the symbol level, the action of pseudo-differential operators on $\cup_k I^k(\Lambda_0\cup\Lambda)$
is by multiplication of symbols, consistent with the standard calculus.
See Theorem 4.13 and formula (5.2) of \cite{MelroseUhlmann79intersection}
for exact sequences which encode the symbol calculus.

\begin{proposition}
\label{prop-MU-Lagr-intsect-calc}
Suppose $P\in\Psiphg^m(X;E)$ is of real principal type with principal symbol $p$,
endowed with a Hamilton field $H=H_q$, and $\tilde p p=q\Id$.
Suppose $\Lambda_0$ is a conic Lagrangian manifold such that $H$ is transversal
to $\Lambda_0$ at $\Lambda_0\cap\Char P$.
Assume that the $H$-flowout of $\Lambda_0\cap\Char P$ defines a conic Lagrangian manifold $\Lambda$
which intersects $\Lambda_0$ in the boundary of $\Lambda$, $\partial\Lambda=\Lambda_0\cap\Char P$.
Let
\[ B\in I^{m+k-1/2}(X,\Lambda_0;\hDB{X}\otimes E) \]
with principal symbol $b$.
There exists $A\in I^k(X,\Lambda_0\cup\Lambda;\hDB{X}\otimes E)$ such that $PA-B\in\Cinfty$.
The principal symbol $a$ of $A$ on $\Lambda\setminus\Lambda_0$ satisfies $pa=0$,
and it is the solution of the transport equation $i^{-1} D_H a + p_s a = 0$ with initial condition
\begin{equation}
\label{eq-transport-initval}
a|\sigma\wedge \intd t|^{-1/2}
  = (2\pi)^{1/4}e^{\pi i/4} \tilde p b|\sigma\wedge \intd q|^{-1/2} (Ht)^{1/2}
\end{equation}
at $\partial\Lambda$, the notation being as in \eqref{eq-iLagr-compat}.
\end{proposition}
\begin{proof}
To construct $A$, we proceed symbolically.
Given, in addition, $C\in I^{m+k-1}(\Lambda_0\cup\Lambda)$, we first solve
\begin{equation}
\label{eq-PA-BC}
PA \equiv B+C\mod I^{m+k-3/2}(\Lambda_0)+I^{m+k-2}(\Lambda_0\cup\Lambda).
\end{equation}
Denote by $c$ the principal symbol of $C$ on $\Lambda\setminus\Lambda_0$.
Let $a$ be the solution of the transport equation
$i^{-1} D_H a + p_s a = \tilde p c$ on $\Lambda$
which satisfies the initial condition \eqref{eq-transport-initval} at $\partial\Lambda$.
Since $p\tilde p=0$ on $\Char P$, $pa=0$ holds at $\partial \Lambda$, hence on all of
$\Lambda$ by Lemma~\ref{lemma-homog-transp-eq}.
Put $a_0=p^{-1}b$ on $\Lambda_0\setminus\Char P$.
The pair $(a_0,a)$ satisfies the compatibility condition \eqref{eq-iLagr-compat}.
Choose $A\in I^k(\Lambda_0\cup\Lambda)$ with principal symbol $(a_0,a)$.
By the symbol calculus, $PA-B\in I^{m+k-1}(\Lambda_0\cup\Lambda)$.
By Theorem~\ref{theorem-amp-transport-eqn} and the properties of $a$,
the $\Lambda$-component of the principal symbol of $PA-B-C$ ranges in $\kernel \tilde p=\image p$.
Therefore, we can find $A'\in I^{k-1}(\Lambda_0\cup\Lambda)$ such that
the $\Lambda$-component of the principal symbol of
$P(A+A')-B-C\in I^{m+k-1}(\Lambda_0\cup\Lambda)$ vanishes.
It follows from the symbol calculus that \eqref{eq-PA-BC} holds with $A$ replaced by $A+A'$.

We use \eqref{eq-PA-BC} to recursively construct $A_{j'}\in I^{k-j'}(\Lambda_0\cup\Lambda)$ such that
\[
P\sum\nolimits _{0\leq j'<j} A_{j'} \equiv B
   \mod I^{m+k-j-1/2}(\Lambda_0)+I^{m+k-j-1}(\Lambda_0\cup\Lambda)\subset I^{m+k-j}(\Lambda_0\cup\Lambda).
\]
Asymptotic summation completes the proof.
\end{proof}

Replacing $X$ by $X\times X$, $\Lambda_0$ by the twisted conormal bundle of the diagonal in $X\times X$,
and $B$ by the identity operator $\Id$, the foregoing arguments lead to a forward parametrix $E_+$.
So, $PE_+\equiv \Id$ holds modulo a smoothing operator,
and $E_+$ is an FIO associated with an intersecting pair of canonical relations.

\section{Boundary parametrices of elastodynamics}
Let $(M,g,C)$ be an elastic body.
Assuming that the operator of elastodynamics, $P=L-\rho D_t^2$, is of real principal type,
we shall construct, microlocally at non-glancing regions, outgoing and incoming Dirichlet parametrices
at the boundary of space-time $X=M\times\R$.
Furthermore, we define the corresponding microlocal DN maps which map displacements to tractions.

The boundary $\partial M$ is endowed with its natural metric tensor and connection.
Let $n$ be the interior unit normal field at $\partial M$.
Discarding a closed subset of the interior of $M$, removes the interface $N$ from consideration,
and $\kappa:(r,y)\mapsto x=\exp_y(r\,n(y))$ maps $[0,R[\times\partial M$ diffeomorphically onto $M$.
We assume that $r$ is the distance from $x$ to $\partial M$,
and that $y=b(x)$ is the unique point in $\partial M$ nearest to $x$.
Parallel transport along boundary orthogonal geodesics, $r\mapsto \kappa(r,y)$,
gives an orthogonal bundle isomorphism between $TM$ and the pullback by $b$ of $T_{\partial M}M$.
Using the isomorphism, we identify sections of the tangent bundle $TM$ with smooth maps
from the interval $[0,R[$ into the space of sections of $T_{\partial M}M$.
Analogously, tensor fields over $M$ are identified with $r$-dependent sections
of the restrictions to $\partial M$ of corresponding tensor bundles.
We equip the bundles with the pullback connections under the inclusion map $\partial M\hookrightarrow M$.
Using the identifications, the elasticity operator $L$ and other differential operators are
regarded as elements of the algebra generated by $D_r=i^{-1}\partial_r$ and the tangential differential operators.
An operator is said to be tangential if it commutes with multiplication by $r$.
We use the complexifications of bundles without distinguishing in notation from the real case.

Covariant derivation in directions tangent to boundary orthogonal geodesics is given by
\[ \nabla_{V} u(x)= \partial_r u(r,y), \quad x=\kappa(r,y), \quad V= \PT_{x\from y}n(y). \]
Put $\nu={\intd r}|_{\partial M}$.
Define $B$ by
\begin{equation}
\label{eq-covar-der-collar}
-i\nabla u = D_r u \otimes \nu + Bu.
\end{equation}
The operator $B$ is tangential.
More precisely, $B$ consists of a family, smoothly parametrized by $r$, of differential
operators which map sections of $T_{\partial M}M$ to sections of $\Hom(T\partial M,T_{\partial M}M)$.
Here $T_{\partial M}M=\R n\oplus T\partial M$ and the corresponding orthogonal decomposition of $TM$ is used.
Notice that $B|_{r=0}=-i\nabla^{T_{\partial M}M}$.

As in \cite[Proposition~11]{SH11rqm}, the elasticity and the traction operator are
\begin{equation}
\label{eq-L-and-T-in-collar}
L = (D_r -i\tr S)(A_0D_r + A_1) + A_1^\sharp D_r + A_2,
\quad
iT  = A_0 D_r +A_1.
\end{equation}
The Hessian $S=\nabla^2 r$ is the shape operator associated with the level surfaces of $r$.
The proof of \eqref{eq-L-and-T-in-collar} uses \eqref{eq-Green-id-elasticity}
with strain tensors replaced by covariant derivatives:
\begin{equation}
\label{eq-Green-id-elast-covar}
\int_{r\geq 0}\int_{\partial M} C\nabla u \overline{\nabla v}\, J\intd V_{\partial M}\intd r
   = \int_{r\geq 0}\int_{\partial M} Lu\,\bar v \, J\intd V_{\partial M}\intd r
      + \int_{\partial M} Tu\,\bar v \, \intd V_{\partial M}.
\end{equation}
Notice that $\kappa^*\intd V_M = J\intd V_{\partial M}\intd r$ holds
with $J$ given by $\partial_r \log J=\tr S$, $J|_{r=0}=1$.
Recall that the bundle identification by parallel transport preserves inner products.
Insert \eqref{eq-covar-der-collar} into \eqref{eq-Green-id-elast-covar}
and perform partial integrations to get \eqref{eq-L-and-T-in-collar}.
Put $B^\sharp=J^{-1}B^* J$ with $B^*$ the adjoint of $B$.
The coefficients in \eqref{eq-L-and-T-in-collar} are
\begin{equation}
\label{eq-Aj-op}
A_0 = \nu\cdot C \cdot\nu, \quad A_1 = \nu\cdot C B, \quad A_1^\sharp = B^\sharp C\cdot\nu, \quad A_2 = B^\sharp C B.
\end{equation}
The dots to the left and right of the stiffness tensor $C$ denote contractions with one of the two left 
and one of the two right slots of $C$, respectively.
The operators $A_j$ and $A_j^\sharp$ are of order $\leq j$ with real principal symbols.
The principal symbols $a_0(r,y)$ and $a_2(r,y,\eta)$, $0\neq\eta\in T_y^*\partial M$,
of $A_0$ and $A_2$ are symmetric positive definite.
The principal symbol of $A_1^\sharp$ is the transpose $a_1(r,y,\eta)^t$ of the principal symbol of $A_1$.
Coordinatewise, at $r=0$, the principal symbols are given by
\begin{equation}
\label{eq-aj-coord-rep}
a_0^{ik} = \nu_j C^{ijkm} \nu_m, \quad a_1^{ik} = \nu_j C^{ijkm} \eta_m, \quad a_2^{ik} = \eta_j C^{ijkm} \eta_m.
\end{equation}

Suppose $P= L -\rho D_t^2$ is of real principal type.
Endow $P$ with the Hamilton field $H$ which satisfies $Ht=1$ on $\Char P$.
Fix
\[ \gamma=(\eta,\tau)\in T_{(y,t)}^*(\partial M\times\R). \]
In general, there are several bicharacteristics which hit the space-time boundary over $\gamma$.
Each corresponds to a unique point in $\Char P\cap (\gamma+\R\nu)$, $\nu=\intd r(y)$.
Suppose $\gamma$ non-glancing, that is, $Hr$ never vanishes at the intersection.
Set $E=\End(T_y M)$, and
\[ a(s)= p(\gamma+s\nu) = a_0 s^2 +\big(a_1(\eta)+ a_1(\eta)^t\big)s+a_2(\eta)-\rho\tau^2 \in\End(E), \]
where $p$ denotes the principal symbol of $P$.
We are going to apply results of Appendix~\ref{sect-spectral-factorization} to
the quadratic self-adjoint polynomial function $a(s)$.
The real spectrum of $a(s)$ equals $\Char P\cap (\gamma+\R\nu)$.
Let $a'(s)$ denote the derivative of $a(s)$.
Define the \emph{outgoing spectrum} (resp.\ \emph{incoming spectrum}) at $\gamma$
as the set $\sigma_\outgoing$ (resp.\ $\sigma_\incoming$) of $s\in\C$, $\IM s\geq 0$,
such that $\kernel a(s)\neq 0$, and $-\tau a'(s)$
is positive (resp.\ negative) definite on $\kernel a(s)$ when $s$ is real.
Proposition~\ref{prop-sqp-factorization}
applies to give spectral factorizations at $\gamma$:
\begin{equation}
\label{eq-out-in-spec-factorization}
a(s)= (s-q_{\outgoing/\incoming}^\sharp) a_0 (s-q_{\outgoing/\incoming}),
\end{equation}
$\spec(q_{\outgoing/\incoming})= \sigma_{\outgoing/\incoming}$, and
$\spec(q_{\outgoing/\incoming}^\sharp)\cap \spec(q_{\outgoing/\incoming})= \emptyset$.
Bicharacteristic strips are outgoing from (resp.\ incoming to) the boundary
iff the distance $r$ increases (resp.\ decreases) as time increases.
The outgoing and incoming spectral factorizations are distinct only if the real spectrum of $a(s)$ is non-empty.
The terms outgoing and incoming spectral factorization are justified by the following lemma.

\begin{lemma}
\label{lemma-p-on-affineline-signs}
Suppose $\gamma+s\nu\in\Char P$.
The restriction of $-\tau a'(s)$ to the null space of $a(s)$
is positive (resp.\ negative) definite iff $Hr(\gamma+s\nu)$ is positive (resp.\ negative).
\end{lemma}
\begin{proof}
Suppose $a(s)w=0$ with $w\neq 0$.
Using the real principal type assumption, we extend $w$ to a neighborhood of $\gamma+s\nu$
as a smooth section which satisfies $pw=0$ on $\Char P$.
The real-valued function $q=(w|pw)$ vanishes on $\Char p$ and is,
in a neighborhood of $\gamma+s\nu$, a defining function of $\Char P$.
Indeed, $\partial/\partial\tau$ is transversal to $\Char P$,
and the Hamilton field $H_q$ of $q$ satisfies
$H_q t=\partial q/\partial\tau =-2\rho\tau (w|w)\neq 0$ on $\Char P$.
Thus $H=cH_q$ holds with $c\tau<0$.
The assertion follows from $Hr=cH_q r=c(w|a'(s) w)$.
\end{proof}

Microlocally at non-glancing points we decompose $P$ into a product of first order operators.
The factorization will be in the algebra generated by $D_r$ and
the algebra of tangential pseudo-differential operators, $\Psitang^\infty=\cup_m \Psitang^m$.
If $y$ denotes local coordinates of $\partial M$, then tangential pseudo-differential
operators are of the form $A(r,y,t,D_y,D_t)$.
We ignore the issue that these operators may not be pseudo-differential near boundary conormals
because we apply the operators only to distributions which are $\Cinfty$-maps
from $r\geq 0$ into distribution sections over the boundary.
Equip the space-time boundary $\partial M\times \R$ with the pullback connection
obtained from the Levi-Civita connection of the boundary $\partial M$ by projecting off the time axis.
We use the geometric symbol calculus for pseudo-differential operators on $\partial M\times \R$.

\begin{proposition}
In a microlocal neighborhood of a given non-glancing boundary point,
there exist $Q, Q^\sharp\in\Psitang^1$ such that
\begin{equation}
\label{eq-eladyn-op-factorized}
L -\rho D_t^2 \equiv (D_r -Q^\sharp)A_0 (D_r - Q)
\end{equation}
holds modulo $\Psitang^{-\infty}+ \Psitang^{-\infty}D_r$.
The principal symbol $q$ of $Q$ and the principal symbol $q^\sharp$ of $Q^\sharp$ have disjoint spectra.
The bicharacteristic strips of $D_r-Q$ are the outgoing bicharacteristics of $L-\rho D_t^2$.
\end{proposition}

The factor $D_r-Q$ is of real principal type because $L-\rho D_t^2$ is.
We call $Q_{\outgoing}=Q$ the outgoing right root of $L-\rho D_t^2$.
Similarly, there exists an incoming right root $Q_{\incoming}$.
\begin{proof}
On the principal symbol level, \eqref{eq-eladyn-op-factorized} necessitates
\begin{equation}
\label{eq-princsymb-eladyn-op-factorized}
a_0 s^2 +\big(a_1(\eta)+ a_1(\eta)^t\big)s+a_2(\eta)-\rho\tau^2
   = \big(s -q^\sharp(\eta,\tau)\big)a_0 \big(s-q(\eta,\tau)\big).
\end{equation}
To find $q$, we use the outgoing factorizations \eqref{eq-out-in-spec-factorization},
that is, we choose $q=q_{\outgoing}$.
The real principal type property and the non-glancing assumption insures
that real eigenvalues $s$ of $a(s)=p(\gamma+s\nu)$
are smooth functions of $\gamma=(y,t,\eta,\tau)$ and of $r\geq 0$, provided $r$ is small.
Furthermore, the sign type of $s$ remains locally constant;
see Lemma~\ref{sqp-residue-semidef} for the positive or negative type of a real eigenvalue.
The integral formula \eqref{sqp-Qr-intformula} implies that the right root $q=q(r,y,\eta,\tau)$
is smooth and is homogeneous of degree $1$ in $(\eta,\tau)$.
The same is true of the left root $q^\sharp$.

Passing from symbols to operators,
\eqref{eq-L-and-T-in-collar} and \eqref{eq-princsymb-eladyn-op-factorized} give
\[ L -\rho D_t^2 = (D_r -Q^\sharp)A_0 (D_r - Q) + R_0D_r+R_1 \]
with $R_j\in\Psitang^j$.
The spectra of $q$ and $q^\sharp$ are disjoint.
Hence the equations $s_jq-q^\sharp s_j=r_j$ have unique solutions $s_j$.
The principal symbol $r_j$ of $R_j$ is homogeneous of degree $j$.
So, $s_j$ is a homogeneous symbol of degree $j-1$.
Choose $S_j\in\Psitang^{j-1}$ with principal symbols $s_j$.
Replacing $Q$ by $Q-A_0^{-1}(S_0Q+S_1)$
and $Q^\sharp$ by $Q^\sharp+(Q^\sharp S_0+S_1)A_0^{-1}$,
the orders of the error terms are decreased by one step, i.e.\ $R_j\in\Psitang^{j-1}$ holds.
Continuing this procedure we find, for every $k\in\N$, $Q$ and $Q^\sharp$ such that $R_j\in\Psitang^{j-k}$ holds.
The proof is completed using asymptotic summation.
\end{proof}

As in \eqref{eq-eladyn-op-factorized}, let $Q=Q_{\outgoing}$ be the outgoing right root
near a conic non-glancing set $\Gamma\subset T^*(\partial M\times \R)$ with compact base.
The principal symbol $q=q(r,y,t,\eta,\tau)$ has a spectral decomposition $q=\sum_s s\psi_s + \psi$.
Here $\psi_s$ projects onto the eigenspace of the outgoing real eigenvalue $s$,
and the spectrum of $\psi$ is contained in the positive upper half-plane.
Following \cite{Taylor75reflection},
we find that $(D_r-Q)W\equiv W(D_r-Q')$ holds with
$W\in\Psitang^0$ elliptic and $Q'\in\Psitang^1$ which has a block diagonal principal symbol with blocks
corresponding to the spectral decomposition of $q$.
Hence, microlocally near the given non-glancing point, we have a parametrix $U$ of the Dirichlet problem
\[ (D_r-Q)U\equiv 0, \quad U|_{r=0}\equiv \Id. \]
Furthermore, $U=\sum_s U_s+V$ where $U_s$ is a (tangential) Fourier integral operator associated
with the canonical relation given by the bicharacteristic strips outgoing from $(0,y,t,s,\eta,\tau)$.
The principal symbol of $U_s$ at $r=0$ equals $\psi_s$.
The Poisson type operator $V$ maps boundary data to sections which are $\Cinfty$ in $r>0$.
Microlocally we construct outgoing and incoming boundary parametrices, $U_\outgoing$ and $U_\incoming$,
such that $u=U_{\outgoing/\incoming}f$ solves $Lu-\rho D_t^2u\equiv 0$ and $u|_{r=0}\equiv f$
microlocally near $\Gamma$ if $\WF(f)\subset\Gamma$.
The wavefront set of $u=U_{\outgoing}f$ (resp.\ $u=U_{\incoming}f$)
is disjoint from incoming (resp.\ outgoing) bicharacteristics which pass over $\Gamma$,
justifying the terminology.
In view of \eqref{eq-L-and-T-in-collar}, in the non-glancing region,
the boundary traction is given as follows:
\[
Tu|_{r=0}  = -i(A_0 D_r u +A_1 u)|_{r=0} \equiv Z_{\outgoing/\incoming} f,
    \quad u=U_{\outgoing/\incoming}f.
\]
Here $Z_{\outgoing/\incoming} = -i(A_0 Q_{\outgoing/\incoming}+A_1)|_{r=0}$ is the outgoing/incoming DN-map
which maps boundary traces of displacements (Dirichlet data) to boundary traces of tractions (Neumann data).
The principal symbol $q= q_{\outgoing/\incoming}$ of $Q_{\outgoing/\incoming}$
is the unique right root of the spectral factorization \eqref{eq-princsymb-eladyn-op-factorized}.
The principal symbol of $Z_{\outgoing/\incoming}$ is
$z_{\outgoing/\incoming}=-i(a_0 q_{\outgoing/\incoming}+a_1)$.
We call $z_{\outgoing/\incoming}$ boundary impedance because it corresponds
to the surface impedance tensor in the physics literature, \cite{LotheBarnett85surfwaveimped}.

The components of $M^\circ\setminus N$ extend as manifolds with smooth boundaries.
So the results of the present section not only give Dirichlet parametrices and impedance operators
for the boundary $\partial M$ but also for both sides $N^\pm$ of the interface $N$.
Thus we have, microlocally at non-glancing regions in $T^*(N^\pm\times\R)$,
Dirichlet parametrices $U^{\pm}_{\outgoing/\incoming}$ and impedances $z^{\pm}_{\outgoing/\incoming}$.

\section{Reflection and transmission}
\label{sect-refl-transm}
We turn to the microlocal study of problem \eqref{eq-eladyn-source-problem}
at the boundary and at the interface.
Given sources $h$ and $h_j$ on the boundary $\partial M\times\R$ and on the interior interface $N\times\R$, respectively, 
we look for solutions of $Pu=0$ which satisfy, microlocally near given non-glancing points,
\[ Tu|_{\partial M\times\R} =  h,\quad u|_{N^+\times\R}-u|_{N^-\times\R}=  h_0,\quad Tu|_{N^+\times\R} + Tu|_{N^-\times\R}= h_1. \]
Recall that the jump of tractions across an interface is the sum of tractions of
both sides because our definition of traction involves the oriented normal.
We assume $u$ to be microlocally outgoing; incoming solutions are handled analogously.
Sources are either given as data or arise as traces of incoming waves.
Using the outgoing Dirichlet parametrices,
we make the ansatz $u=U_{\outgoing}f$ and $u=U_{\outgoing}^\pm f^\pm$
at the boundary and at the interface, respectively.
This leads, microlocally in non-glancing regions,
to the equations $Z_{\outgoing}f=h$ on $\partial M\times\R$ and
\[ f^+-f^- = h_0, \quad Z^+_{\outgoing} f^+ + Z^-_{\outgoing} f^- = h_1 \quad\text{on $N\times\R$.} \]
The equation on $N\times\R$ reduces to
$(Z^+_{\outgoing} + Z^-_{\outgoing}) f^+ = h_1 + Z^-_{\outgoing} h_0$.
Therefore, we ask for properties of $Z_{\outgoing}$ and of $Z^+_\outgoing +Z^-_\outgoing$,
in particular, whether these are elliptic operators.

The projection $T^*(\partial M\times\R)\to \partial M$ off fiber and time axis induces,
by pullback of the complexification of the bundle $T_{\partial M}M\to\partial M$, a rank $3$ bundle $E$.
Similarly, we have rank $3$ bundles $E^\pm\to T^*(N\times\R)$
with respect to the sides $N^\pm$ of the interface $N$.
Points in $T^*(N\times\R)\setminus 0$ are said to be non-glancing
if they are non-glancing for both sides of $N$.
Non-glancing sets are open.
Over the non-glancing sets $E$ and $E^\pm$ split into subbundles of locally constant rank:
$E = E_c\oplus E_r$ and $E^\pm = E_c^\pm\oplus E_r^\pm$.
The subbundles with subscripts ${}_c$ and ${}_r$ correspond respectively to the non-real and the real
spectrum of the outgoing spectral factorization;
see \eqref{eq-H-spec-decomp} for the fiberwise decomposition where $Q=q_{\outgoing}$.
The non-glancing set decomposes into three open subsets, the hyperbolic, mixed, and elliptic regions.
The hyperbolic region $\Hyp$ is defined by $E_c=0$ and $E_c^+\cap E_c^- =0$,
and the elliptic region $\Ell$ by $E_c=E$ and $E_c^+=E_c^-=E$.
The complement of $\Hyp\cup\Ell$ in the non-glancing set constitutes the mixed region $\Mix$.
Notice that, at $N$, $\Hyp$ includes any region which is hyperbolic with respect to at least one side of $N$,
e.g., hyperbolic--mixed or hyperbolic--elliptic regions in the sense of \cite{StefUhlVasy21transm} are contained in $\Hyp$.

In this section we study outgoing impedances $z$ and $z^++z^-$ only at $\Hyp\cup\Mix$.
The elliptic region $\Ell$ is analyzed in Section~\ref{sect-surf-waves}.

\begin{proposition}
\label{prop-z-ell-at-Hyp}
$Z_\outgoing$ and $Z_{\outgoing}^+ + Z_{\outgoing}^-$ are elliptic in $\Hyp$.
\end{proposition}
\begin{proof}
It follows from Proposition~\ref{prop-Z-nonsing-if-sigma-real} and
Lemma~\ref{lemma-p-on-affineline-signs} that
$-\tau\IM (u|z_{\outgoing} u) \geq 0$ holds for $u\in E$ with strict inequality if the
$E_r$-component of $u$ is non-zero.
This implies
\[ -\tau\IM \big(u|(z_{\outgoing}^++z_{\outgoing}^-) u\big) > 0 \quad\text{if $u\not\in E_c^+ \cap E_c^-$.} \]
Hence $\kernel z_{\outgoing}\subset E_c$
and $\kernel (z_{\outgoing}^+ + z_{\outgoing}^-)\subset E_c^+\cap E_c^-$.
\end{proof}
No assumption about the symmetry group of the elastic medium is made in Proposition~\ref{prop-z-ell-at-Hyp}.
For an analogous result in $\Mix$, we shall assume the medium isotropic and use
a matrix representation on impedance.

Matrix representations of the impedance $z=z_\outgoing(\eta,\tau)$
are obtained from eigenvectors of the polynomial $a(s)=p(\gamma+s\nu)$:
\begin{equation}
\label{eq-z-on-eigenvec}
izv=(s a_0+a_1)v
\quad\text{whenever $a(s)v=0$ and $s\in\sigma_\outgoing$ hold.}
\end{equation}
Here $q=q_\outgoing(\eta,\tau)$ is outgoing right root introduced in \eqref{eq-out-in-spec-factorization}.
Equivalently, the six-dimensional vector $[v,izv]^t$ is an eigenvector
of the Stroh matrix \eqref{factor-Stroh-matrix}.
In case of non-semisimple eigenvalues, principal vectors must be used in the obvious way.

Suppose the elastic material isotropic with Lam\'e parameters $\lambda, \mu$ 
satisfying $\mu>0$ and $\lambda+\mu>0$.
In view of \eqref{eq-p-iso-eladyn}, we consider the polynomial
\[ a(s)= p(\xi,\tau) = (\lambda+\mu)\xi\otimes\xi +\mu\xi^2\Id -\rho\tau^2\Id, \quad \xi=\eta+s\nu.  \]
Recall that $\eta,\nu\in T_y^* M$ are orthogonal, and $|\nu|=1$.
By \eqref{eq-aj-coord-rep},
\[
a_0= (\lambda+\mu)\nu\otimes\nu + \mu\Id, \quad
a_1(\eta)= \lambda\nu\otimes\eta + \mu\eta\otimes\nu.
\]
The inequalities $\rho\tau^2>(\lambda+2\mu)\eta^2$,
$(\lambda+2\mu)\eta^2>\rho\tau^2>\mu\eta^2$,
and $\mu\eta^2>\rho\tau^2>0$ 
define the hyperbolic, mixed, and elliptic regions, respectively.
Suppose $\tau<0$.
The outgoing eigenvalues $s=s_{p/s}$, $\kernel a(s)\neq 0$, are the roots of $s^2=\tau^2/c_{p/s}^2-\eta^2\neq 0$
which are positive real or positive imaginary.
Choose $\zeta\neq 0$ orthogonal to $\nu$ and $\eta$.
Fix $\tilde s$ so that $\eta+\tilde s\nu$ and $\eta+s_s\nu$ are orthogonal,
i.e., $\eta^2+\tilde s s_s=0$.
A straightforward computation based on \eqref{eq-z-on-eigenvec} gives
\begin{align*}
	\mu^{-1} iz \zeta &= s_s\zeta, \\
	\mu^{-1} iz (\eta+\tilde s\nu)  &= -2\eta^2 \nu +(s_s+\tilde s)\eta, \\
	\mu^{-1} iz (\eta+s_p\nu) &= (s_s^2-\eta^2)\nu +2s_p\eta.
\end{align*}
The formulas represent $z=z_\outgoing$ in a basis of SH-SV-P polarizations.
Because of $s_p\neq \tilde s$, it is easy to obtain the matrix representation of $z$
with respect to the basis $\zeta,\nu,\eta$.

\begin{proposition}
\label{prop-z-ell-at-Mix}
Suppose the elastic body is isotropic.
Then $Z_\outgoing$ and $Z_{\outgoing}^+ + Z_{\outgoing}^-$ are elliptic in $\Mix$.
\end{proposition}
\begin{proof}
Suppose $(\eta,\tau)\in\Mix$.
By the proof of Proposition~\ref{prop-z-ell-at-Hyp} it suffices to show
that $\kernel z \cap E_c =0$ and $\kernel(z^+ + z^-)\cap E_c^+\cap E_c^-=0$.
Since $s_s$ is real, $E_r$ contains the orthogonal complement of $\xi=\eta+s_s\nu$.
Therefore, $\dim E_c\leq 1$ and $\min(\dim E_c^+,\dim E_c^-)\leq 1$.
We consider only $z^++z^-$.
Suppose $\dim E_c^+=1$ and $E_c^+\cap E_c^-\neq 0$.
Thus $s_p^+\in i\R_+$ and $E_c^+\cap E_c^-=\C(\eta+s_p^+\nu)$.
We claim that the following holds:
\[ (z^++z^-)(\eta+s_p^+\nu)\neq 0. \]
Suppose $E_c^-=E$.
Then $0<s_s^-<s_p^-$ in $i\R$.
Define $\hat s$ and $t$ by
\[ \eta^2+\hat s s_s^-=0\quad\text{and}\quad t s_p^-+(1-t)\hat s=s_p^+. \]
Observe that $\hat s\in i\R_+$ and $|s_p^\pm|<|\eta|<|\hat s|$, hence $t>0$.
The inner product of $\eta$ and
\[ i(z^++z^-)(\eta+s_p^+\nu) = iz^+(\eta+s_p^+\nu)+tiz^-(\eta+s_p^-\nu) + (1-t)iz^-(\eta+\hat s\nu) \]
equals $\eta^2$ times
$2s_p^+\mu^+ + \big(t(s_p^--s_s^-)+s_p^++s_s^-\big)\mu^-$
which is positive imaginary.
If $E_c^-\neq E$, then $E_c^-=E_c^+$, and the result holds with $t=1$.
This proves the claim and the lemma.
\end{proof}

\section{Surface waves}
\label{sect-surf-waves}
The elliptic region is the open subset $\Ell\subset T^*(\partial M\times\R)\cup T^*(N\times\R)$
over which no incoming or outgoing bicharacteristics pass.
However, elastic waves satisfying zero traction boundary conditions may have
singularities in $\Ell$ because $\Char Z\cap\Ell$ may be non-empty.
This happens for isotropic media and is the reason for the existence of Rayleigh surface waves.
The free surface wave theory of \cite{LotheBarnett85surfwaveimped} implies,
when translated into the language of microlocal analysis as done below,
that $Z$ is, for any elastic symmetry, of real principal type in $\Ell$.
The elliptic region $\Ell$ is called the subsonic region in Barnett--Lothe theory,
and the part of the glancing region which bounds $\Ell$ is called transsonic.

Suppose $0\neq\eta\in T_y^*(\partial M)$.
Set $\nu=\intd r(y)$ as before.
For real $\tau$, the condition $\gamma=(\eta,\tau)\in \Ell$ is equivalent to the positive definiteness of
\[ p(\gamma+s\nu) = a_0 s^2 +\big(a_1(\eta)+ a_1(\eta)^t\big)s+a_2(\eta)-\rho\tau^2, \quad s\in\R. \]
Decreasing $\tau^2$ does not lead out of $\Ell$.
Put $\tau_\eta=\sup_{(\eta,\tau)\in\Ell}\tau$.
Then $|\tau|=\tau_\eta>0$ describes the part of glancing set which bounds $\Ell$.
Denote by $q=q(\eta,\tau)$ the unique spectral right root of the polynomial $p(\gamma+s\nu)$
with spectrum contained in the upper complex half-plane.
The principal symbol of $Z$ is the impedance $z$
given by $iz(\eta,\tau)=a_0q(\eta,\tau)+a_1(\eta)$.

The following fact is well-known.
\begin{lemma}
\label{lemma-z-posdef}
At $\tau=0$ the impedance $z$ is positive definite.
\end{lemma}
\begin{proof}
Suppose $\tau=0$.
Let $u_0\in E$.
Put $u(r)=e^{irq}u_0$ and $D=-i\partial_r$.
As $r\to\infty$, $u(r)$ decreases to zero exponentially.
Using \eqref{eq-Z-int-over-r} and \eqref{eq-aj-coord-rep}, we get
\begin{align*}
(zu_0|u_0) &= \int_0^\infty (a_0Du|Du)+(a_1(\eta)u|Du)+(Du|a_1(\eta)u)+(a_2(\eta)u|u)\intd r \\
           &= \int_0^\infty (Cw(r)| w(r)) \intd r.
\end{align*}
Here $w(r)=Du(r) \nu+u(r) \eta$ is the symmetrization of the tensor $Du(r)\otimes \nu+u(r)\otimes\eta$.
Since the stiffness tensor $C$ is positive definite, $(zu_0|u_0)\geq 0$ holds.
Suppose $(zu_0|u_0)\leq 0$.
This implies $w=0$.
Since $\nu$ and $\eta$ are orthogonal, $0=(w(r)|\nu\nu)=(Du(r)|\nu)|\nu|^2$ holds.
We infer that $(u(r)|\nu)$ is constant, hence zero.
Suppose $\mu$ is orthogonal to $\nu$.
It follows that $(Du(r)|\mu)|\nu|^2= (w(r)|\nu \mu)=0$.
Thus $Du=0$, and $u_0=0$.
\end{proof}

The boundary value $v=u|_{\partial M\times\R}$ of an elastic wave $u$ which satisfies
the free surface boundary condition, $Tu|_{\partial M\times\R}=0$,
is a solution of the equation $Zv=0$.
Next, we show that the results of Section~\ref{section-Lagr-solns} apply to $Zv=0$ in $\Ell$.

\begin{proposition}
\label{prop-Z-realprinctype}
$Z$ is of real principal type in $\Ell$.
The leading symbol of $Z$ is $z+\hh z_1$,
where $z=-i(a_0q+a_1)$ is the impedance, and $z_1$ is the unique solution of
\begin{equation}
\label{eq-z-lead-symb}
q^*z_1 -z_1q + (\vnabla\cdot\hnabla a_1^t)q - \vnabla q^*a_0\cdot\hnabla q +ia_{2,1}+ i\tr(S)z+i\partial_r z =0.
\end{equation}
Here $S=\nabla^2 r$, $\partial_r z=(\partial z/\partial r)|_{r=0}$,
and $a_2+\hh a_{2,1}$ is the leading symbol of $A_2$.
\end{proposition}
We used $z=-i(a_0q+a_1)$ to define the derivative $\partial_r z$.

\begin{proof}
In $\Ell$, $z$ is self-adjoint.
By Proposition~\ref{prop-Re-Z-posdef}, $\RE z$ is positive definite;
consequently, at least two eigenvalues of $z$ are positive.
Proposition~\ref{prop-Zdot-posdef} implies that $\partial z/\partial \tau^2$ is negative definite.
Using perturbation theory of eigenvalues, e.g.\ \cite[Theorem \RN 2-5.4]{Kato76perturb},
we find that $\partial \det z/\partial\tau\neq 0$ holds at $\Char Z$.
In view of Lemma~\ref{lemma-z-posdef},
either $z(\eta,\tau)$ is positive definite for $|\tau|<\tau_\eta$,
or there exists a positive smooth function $\hat\tau$ on an open neigborhood $W\subset T^*(\partial M)$ of $\eta$,
such that, for $\zeta\in W$, $(\zeta,\tau)\in\Char Z$ holds iff $|\tau|=\hat\tau(\zeta)$.
Clearly, $\dim\kernel z(\eta,\tau)= 1$ holds for $(\eta,\tau)\in\Char Z$.
The real principal type property of $Z$ now follows.

For the leading symbol, we procede as in \cite[Lemma 18]{SH11rqm}.
Recall \eqref{eq-L-and-T-in-collar}.
Comparing coefficients of powers of $D_r$,
we see that \eqref{eq-eladyn-op-factorized} is equivalent to the following two equations between tangential operators:
\begin{align*}
A_0Q +A_1 + Q^\sharp A_0 + A_1^\sharp -i\tr(S)A_0 &= 0, \\
[D_r,A_0Q+A_1] -Q^\sharp A_0 Q +A_2-\rho D_t^2 -i\tr(S) A_1 &=0.
\end{align*}
Multiply the first equation from the right by $Q$, and add the equations:
\begin{equation}
\label{eq-to-determine-leadsymb-of-Z}
A_0Q^2+A_1Q+ A_1^\sharp Q +A_2-\rho D_t^2 = -\tr(S) Z -[D_r,A_0Q+A_1].
\end{equation}
The principal symbol of the first order operator on the right-hand side is $-\tr(S)z-\partial_r z$.
Using Propositions~\eqref{prop-geosymb-composition} and \eqref{prop-geosymb-adjoint},
we compute the leading symbol of the second order operator on the left.
First, we determine the symbols of the operators $A_1$ and $A_1^\sharp$ defined in \eqref{eq-Aj-op}.
The horizontal derivative of the principal symbol of $B|_{r=0}$ vanishes,
because this operator is a covariant derivative.
It follows that the leading symbol of $A_1$ equals its principal symbol $a_1$.
At $r=0$, $A_1^\sharp=A_1^*$ holds.
The leading symbol of $A_1^*$ is $a_1^t -i\hh  \vnabla \cdot \hnabla a_1^t$,
and those of $iZ+A_1^*$ and $Q=A_0^{-1}(iZ-A_1)$ are
\[ -q^*a_0+i\hh (z_1-\vnabla\cdot\hnabla a_1^t) \quad\text{and}\quad q+i\hh a_0^{-1}z_1, \]
respectively.
Therefore, the leading symbol of $A_0Q^2+(A_1+ A_1^*)Q = (iZ+A_1^*)Q$ is 
\[ (a_0q^2+a_1q+a_1^tq)+i\hh\big((z_1-\vnabla\cdot\hnabla a_1^t)q-q^* z_1 +\vnabla q^* a_0\cdot\hnabla q \big).\]
The principal symbol of the second order operator on the left-hand side
of \eqref{eq-to-determine-leadsymb-of-Z} must be zero.
This is confirmed by the solvency equation \eqref{sqp-solvency}.
Equating the lower order term of the leading symbol on the left with the
principal symbol on the right of \eqref{eq-to-determine-leadsymb-of-Z} proves \eqref{eq-z-lead-symb}.
\end{proof}

It is clear from the proof that, for $(\eta,\tau)\in\Ell\setminus\Char(Z)$, 
either $z(\eta,\tau)$ is positive definite or $z(\eta,\tau)$ has exactly one negative eigenvalue.
In the latter case, there exists a unique $0<\tau_R<|\tau|$ satisfying $(\eta,\tau_R)\in\Char(Z)$.
This is the basic criterion of Barnett--Lothe theory for the existence of Rayleigh--type surface waves.
Other useful existence criteria are stated in terms of the limit of $z(\eta,\tau)$ as $\tau\to \tau_\eta$;
see \cite{LotheBarnett85surfwaveimped}.

Surface waves $u$ are, microlocally near $\Ell\cap\Char(Z)$, given by $u=Uv$, where $U$
is the Dirichlet parametrix, and $v$ solves the real principal type system $Zv\equiv 0$ on $\Ell$.

\begin{remark}
At the interior interface $N\times\R$, surface waves of Stoneley--type may occur.
These do arise when $\Char(Z^++Z^-)\cap\Ell\cap T^*(N\times\R)\neq\emptyset$.
In fact, the proof of Proposition~\ref{prop-Z-realprinctype} applies also to $Z^++Z^-$ instead of $Z$.
To see this, notice the positive definiteness of $z^++z^-$ at $\tau=0$,
of $\RE(z^++z^-)$ and of $-\partial(z^++z^-)/\partial \tau^2$ in $\Ell\cap T^*(N\times\R)$.
\end{remark}

\section{Propagation of polarization}
\label{section-propag-polar}

Suppose $P=L-\rho D_t^2$ of real principal type.
Endow $\Char P$ with the Hamilton field $H$ determined by $Ht=1$.
We say that a Lipschitz continuous map $\gamma:I\to\bTstar(M\times\R)$ from an interval $I\subset\R$
into the $b$-cotangent bundle is an outgoing broken bicharacteristic of \eqref{eq-eladyn-source-problem}
iff $\gamma(I)$ does not intersect the glancing set and, at time $t\in I$, one of the following holds:
\begin{compactenum}[(i)]
\item $\gamma'(t)=H(\gamma(t))$, or
\item\label{item-bichar-refl-transm}
$\gamma(s) \in T^*(M'\times\R)\cap\Char P$ for $s-t>0$ small.
\end{compactenum}
Case \eqref{item-bichar-refl-transm} arises at points of reflection or transmission.
We call the sequence of components of $\Char P$ containing bicharacteristic
segments the signature of a broken bicharacteristic $\gamma$.

We look for outgoing solutions $u$ of \eqref{eq-eladyn-source-problem} which are generated by Lagrangian sources.
By outgoing we mean that time $t$ is bounded from below on $\supp(u)$.
Since the boundary and the interface are non-characteristic,
\[ \bWF(u)\subset T^*(M'\times\R)\cup T^*(\partial M\times\R)\cup T^*(N\times\R)\subset \bTstar(M\times\R) \]
has to hold.
Here we use the abbreviation $M'=M\setminus (\partial M\cup N)$ from Section~\ref{sect-linear-elast}.

By Proposition~\ref{prop-z-ell-at-Mix}, the first assumption of the following proposition
is true in the case of isotropic elastic materials.

\begin{proposition}
Suppose $Z_\outgoing$ and $Z^+_\outgoing + Z^-_\outgoing$ are elliptic in $\Mix$.
Suppose $f\in I^{k+3/2}(\Lambda_f)$ is a compactly supported Lagrangian distribution in $M'\times\R$,
and $H$ is transversal to the Lagrangian manifold $\Lambda_f$ at $\Lambda_f\cap\Char P$.
Suppose every broken bicharacteristic issuing from $\Lambda_f\cap\Char P$ extends until a time $t_1$.
There is a unique outgoing solution $u(x,t)$, $t<t_1$,
of the transmission problem \eqref{eq-eladyn-source-problem}
which satisfies homogeneous boundary and transmission conditions, $h=h_j=0$.
Moreover, $u$ is a Lagrangian distribution with associated Lagrangian manifolds
labelled by signatures of the broken bicharacteristics issuing from $\Lambda_f$.
Along any broken bicharacteristic the polarization is completely determined by 
transport equations \eqref{eq-transport-polariz} along segments, and by
initial values at $\Lambda_f\cap\Char P$ and at points on the boundary or the interface
in terms of reflection and transmission laws; see \eqref{eq-reflection-law} below.
\end{proposition}
\begin{proof}
Applying Proposition~\ref{prop-MU-Lagr-intsect-calc}, we obtain
$w\in I^{k}(\Lambda_f\cup\Lambda)$ such that $Pw\equiv f$ holds modulo a function
which is $\Cinfty$ except for jump discontinuities across the interface $N\times\R$.
Restriction to a hypersurface is an FIO of order $1/4$.
The Cauchy data of $w$, that is the restrictions of $w$ and the traction $Tw$ to the boundary and the interface,
are Lagrangian distributions of orders $k+1/4$ and $k+5/4$.
The Cauchy data of $w$ are related by the incoming DN operators $Z_\incoming$.
Note that, because $\Lambda$ does not intersect the glancing set,
the Lagrangian manifold associated with the Cauchy data
is the restriction of $\Lambda$ to $\Hyp\cup\Mix$.
Replace $u$ by $u+w$.
Now we have a transmission problem \eqref{eq-eladyn-source-problem}
with $f\equiv 0$, and $h,h_j$ are, in $\Hyp\cup\Mix$, given by the Cauchy data of $w$.
Moreover, $h$ and $h_1$ are Lagrangian distributions of order $k+5/4$, and $h_0$ of order $k+1/4$.
Following the outline in Section~\ref{sect-refl-transm}, we set
\begin{equation}
\label{eq-u-of-sources-h-hj}
u_1 = UZ^{-1}h + U^+(Z^++Z^-)^{-1}(h_1+Z^- h_0) + U^-(Z^++Z^-)^{-1}(h_1-Z^+ h_0).
\end{equation}
Here, inverses are parametrices of the elliptic operators $Z$ and $Z^++Z^-$,
and the operators $U, U^\pm$ are outgoing Dirichlet parametrices.
The canonical relations of the FIO parts of $U$ and $U^\pm$ are extended by Hamilton flowout until 
the boundary or the interface are hit again or the final time $t_1$ is reached.
The Lagrangian manifold generated by the flowouts is $\tilde\Lambda$.
Now $u=u_1\in I^{k}(\tilde \Lambda)$ solves \eqref{eq-eladyn-source-problem} for $f\equiv 0$
and modified boundary sources $h,h_j$.
The infimum of $t$ over
\( (\WF h\cup \WF h_0\cup\WF h_1)\cap(\Hyp\cup\Mix) \)
has strictly increased after modification.
The number of reflections and transmissions along outgoing broken
bicharacteristics issuing from the wavefront sets of the sources is bounded.
Continuing as above, we reach, after finitely many steps, time $t=t_1$.
Thus we have shown that, under our assumptions, \eqref{eq-eladyn-source-problem}
can be solved modulo $\Cinfty$ errors for Lagrangian sources.
By the results of Section~\ref{sect-linear-elast} we obtain exact solutions.
\end{proof}

Let $q_{\outgoing}^\pm=\psi^\pm+\sum\nolimits _s s\psi_s^\pm$ be the spectral decompositions of outgoing right roots,
where $\psi_s^\pm$ is the projector onto the eigenspace with real eigenvalue $s$
and $\psi^\pm$ maps into $E_c^\pm$.
Apply the symbol calculus to the construction \eqref{eq-u-of-sources-h-hj}.
For incoming polarization over $N^+$,
the initial values of reflected and transmitted polarizations at the interface 
can be concisely written as
\begin{equation}
\label{eq-reflection-law}
a^{\pm}_{\outgoing,s}
  = -\psi_{s}^{\pm} (z^+_\outgoing+z^-_\outgoing)^{-1} (z^+_\incoming \pm z^{\mp}_\outgoing)a^+_\incoming.
\end{equation}
An analoguous but simpler law holds at the boundary.

\begin{remark}
Our assumption that $Z$ be elliptic in the mixed region means that we do not cover the
phenomenon of supersonic surface waves
in special tranversely isotropic media, \cite{GundWangLothe91secluded}.
\end{remark}

Finally we turn to surface waves.
In the elliptic region, $Z$ and $Z^++ Z^-$ are of real principal type, possibly even elliptic.
Endow these operators with a Hamilton field, also denoted $H$, which satisfies $Ht=1$.
\begin{proposition}
\label{prop-surface-wave-propagation}
Let $h\in I^{k+1/2}(\Lambda_h)$ a compactly supported Lagrangian distribution on $\partial M\times \R$.
Suppose that $\Lambda_h\subset\Ell\cap T^*(\partial M\times\R)$ and $H$ is transversal
to $\Lambda_h$ at $\Lambda_h\cap\Char Z$.
Suppose that no (Rayleigh) bicharacteristic issuing from $\WF(h)$ intersects the glancing set before time $t_1$.
The transmission problem \eqref{eq-eladyn-source-problem} with $f=0=h_j$ has a unique
outgoing solution $u(x,t)$, $t<t_1$.
Furthermore, $u\equiv Uw$ modulo $\Cinfty$, where $U$ is the microlocal Dirichlet
parametrix at $\Ell\cap T^*(\partial M\times\R)$, and $w\in\Dprime(\partial M\times\R)$ is the microlocally
outgoing solution of $Zw\equiv h$.
The polarization of $w$ satisfies a transport equation \eqref{eq-transport-polariz}
with initial values determined from the principal symbol of $h$ at $\Lambda_h\cap\Char Z$.
\end{proposition}
\begin{proof}
We apply Proposition~\ref{prop-MU-Lagr-intsect-calc}.
We obtain
$w\in I^{k}(\Lambda_h\cup\Lambda)$, $\Lambda\subset\Char Z$ the Hamilton flowout of
$\Lambda_h\cap\Char Z$, such that $Zw\equiv h$ holds for $t<t_1$.
Observe that $Uw$ satisfies the transmission problem modulo smooth error.
Using Section~\ref{sect-linear-elast}, we obtain an exact solution $u\equiv Uw$.
\end{proof}

Of course, actual surface waves can only arise when $\Lambda_h\cap\Char Z$ is non-empty.
To set up the transport equation for the polarization of $u|_{\partial M\times\R}\equiv w$
it is necessary to evaluate the subprincipal symbol of $Z$ using
the formulas \eqref{eq-subprinc-symb} and \eqref{eq-z-lead-symb} with $p=z$ and $p_1=z_1$.
This is cumbersome but algorthmically straightforward.
Observe that the subprincipal symbol of $Z$, hence also the polarization of $w$,
depends not only on the elasticities on $\partial M$
but also on their $r$-derivatives in direction transversal to the boundary.
Furthermore, the polarization depends on the curvature of the boundary.

In elasticity, dispersion of Rayleigh wave velocity has been observed,
i.e., velocity depends on frequency and $r$; compare \cite{MaNaTaWa15dispersion}.
Dispersion refers to $u$ and not its boundary trace $w$, because the latter
is Lagrangian with frequency-independent wave speed.

If $\Char(Z^++Z^-)\cap\Ell\neq\emptyset$, then
existence of Stoneley type waves at the interior interface is proved
in the same way as in Proposition~\ref{prop-surface-wave-propagation}.
See the recent analysis \cite{Zhang2020rayleigh} of Rayleigh and Stoneley waves in isotropic media.

\appendix
\section{Elastic symmetries}
\label{sect-elastic-symmetries}
Elastic materials are characterized by the symmetry group of the stiffness tensor.
In this section we recall elastic symmetries.
For transversely isotropic media which are small perturbations from isotropy,
we prove a result about the eigenvalues of the acoustic tensor.

Fix a point in a $3$-dimensional elastic body.
With respect to an orthonormal basis of the tangent space at that point,
the components of the stiffness tensor $C$ are the elasticity constants $C^{ijkm}$.
As in Section~\ref{sect-linear-elast} we view $C$ as a contravariant $4$-tensor.
The symmetry properties of $C$ can be expressed as $C^{ijkm} = C^{jikm}=C^{km ij}$.
The dimension of the vector space $E$ of stiffness tensors equals $21$.

The special orthogonal group $SO(3)$ acts on tensors.
More specifically, $O\in SO(3)$ acts $2$-tensors, viewed as $3\times 3$-matrices, by $O\cdot S=OSO^t$,
and by $(O\cdot C)[S] = O^t C[OSO^t] O$ on $C\in E$.
From the representation theory of $SO(3)$ the irreducible (real) representations are known.
These are isomorphic to left translation on the spaces $H_n$ of harmonic polynomial
functions homogeneous of degree $n$.
The dimension of $H_n$ is $2n+1$.
Note that $H_2$ can be identified with the space of traceless symmetric matrices,
and $H_4$ with the space of the symmetric $4$-tensors $C$ which satisfy $C^{iikm}=0$ for all $k$ and $m$.
The representation on stiffness tensors decomposes into a direct sum of irreducibles:
$E \simeq H_4\oplus H_2\oplus H_2 \oplus H_0\oplus H_0$.
More explicitly,
\begin{equation}
\label{eq-c-HABlambdamu}
\begin{aligned}
C^{ijkm} &= H^{ijkm} + \big(\delta^{ij}A^{km}+\delta^{km}A^{ij}\big)
   +\big(\delta^{ik}B^{jm} + \delta^{jm}B^{ik} + \delta^{im}B^{jk} + \delta^{jk}B^{im} \big) \\
   &\phantom{==} +\lambda\delta^{ij}\delta^{km} +\mu\big(\delta^{ik}\delta^{jm}+\delta^{im}\delta^{jk}\big).
\end{aligned}
\end{equation}
Here $H\in H_4$, $A,B\in H_2$, and $\lambda,\mu\in H_0=\R$.
Observe that
\begin{align*}
C^{iikk} &= 9\lambda + 6\mu, \quad C^{ijij} = 3\lambda + 12\mu, \\
C^{iikm} &= 3A^{km} + 4 B^{km} + (3\lambda +2\mu)\delta^{km}, \\
C^{ijim} &= 2A^{jm} + 5 B^{jm} + (\lambda +4\mu)\delta^{jm},
\end{align*}
which implies that $C$ determines $\lambda,\mu, A, B, H$ uniquely.

For $C\in E$, the stabilizer $G(C)\subset SO(3)$ consists of the matrices $O$ which satisfy $O\cdot C=C$.
Each stabilizer is conjugate to one of the following $8$ subgroups:
$SO(3), O(2)$, the octahedral group $\mathcal{O}$, the dihedral groups $D_4, D_3, D_2$,
the cyclic group $Z_2$, or the trivial subgroup.
See \cite{ForteVian96symela} for details.

Suppose that $G(C)$ contains a subgroup conjugate to $SO(2)$ with axis of rotation a unit vector $J$.
Then the elastic material is either isotropic, $G(C)=SO(3)$,
or transversely isotropic with respect to $J$, i.e.\ $G(C)$ is conjugate to $O(2)$.
The tensors $A,B,H$ in \eqref{eq-c-HABlambdamu} are invariant under rotations around the axis $J$.
This implies that $A$ and $B$ are scalar multiples of $J\otimes J-1/3$,
and $H$ equals $J\otimes J\otimes J\otimes J$ modulo terms which can be absorbed into $A,B,\lambda,\mu$.
So \eqref{eq-c-HABlambdamu} becomes
\begin{equation}
\label{eq-C-transverse-iso}
\begin{aligned}
C^{ijkm} &= \lambda\delta^{ij}\delta^{km} +\mu\big(\delta^{ik}\delta^{jm}+\delta^{im}\delta^{jk}\big)
   + \alpha\big(\delta^{ij}J_kJ_m+\delta^{km}J_iJ_j\big) \\
  &\phantom{==}
   +\beta \big(\delta^{ik}J_jJ_m + \delta^{jm}J_iJ_k + \delta^{im}J_jJ_k + \delta^{jk}J_iJ_m \big)
   +\gamma J_iJ_jJ_kJ_m.
\end{aligned}
\end{equation}
Compare \cite[\S 6]{Synge57elast}.
Note that $C$ is determined by five real parameters $\lambda,\mu,\alpha,\beta,\gamma$ and by the axis of rotation, $J$.
When $J=(0,0,1)$ the elasticities
$C^{1111}$, $C^{1122}$, $C^{1133}$, $C^{2323}$, $C^{3333}$
can be used as parameters; see \cite{Tanuma07Stroh}.
In the isotropic case, only the Lam\'e parameters $\lambda$ and $\mu$ can be non-zero.
Positive definiteness (strong convexity) of an isotropic stiffness $C$ holds iff $\mu>0$ and $\lambda+\mu>0$.

Given a stiffness tensor $C$ and a (co-)vector $\xi$,
set $\ell(\xi)^{ik}=C^{ijkm}\xi_j\xi_m$.
The tensor $\ell(\xi)$, called acoustic tensor, is a symmetric endomorphism of the tangent space.
Suppose $C$ is isotropic.
Then, writing $\xi^2=|\xi|^2$,
\begin{equation}
\label{eladyn-iso-acoustic-tensor}
\ell(\xi) = \ell_{\lambda,\mu}(\xi)
   = (\lambda+\mu)\xi\otimes\xi +\mu\xi^2\Id = \xi^2\big((\lambda+2\mu)\pi_p(\xi) +\mu \pi_s(\xi)\big).
\end{equation}
Here $\pi_p(\xi)$ and $\pi_s(\xi)$ are, for $\xi\neq 0$,
the orthoprojectors onto the line spanned by $\xi$ and onto its orthogonal complement, respectively.

Suppose $C$ is transversely isotropic.
Using \eqref{eq-C-transverse-iso}, the acoustic tensor is found to be
\begin{equation}
\label{eladyn-transiso-acoustic-tensor}
\begin{aligned}
\ell(\xi) &= \ell_{\lambda,\mu}(\xi) + t(\xi), \\
t(\xi) &=(\alpha+\beta)(J|\xi)(\xi\otimes J+J\otimes\xi) +\beta(J|\xi)^2I+(\beta \xi^2+\gamma(J|\xi)^2)J\otimes J.
\end{aligned}
\end{equation}
The subspace orthogonal to $\xi$ and $J$ is an eigenspace of $\ell(\xi)$
with eigenvalue $\mu \xi^2+\beta(J|\xi)^2$, \cite[(6.7)]{Synge57elast}.
For small parameters $\alpha,\beta,\gamma$, we regard \eqref{eladyn-transiso-acoustic-tensor} 
as a perturbation of an isotropic elastic system.
\begin{proposition}
\label{prop-transv-iso-perturb}
Suppose $\lambda +\mu\neq 0$, $\pi_s(\xi)J\neq 0$, and $\beta+\gamma(J|\hat \xi)^2\neq 0$ where $\hat \xi=\xi/|\xi|$.
The eigenvalues of \eqref{eladyn-transiso-acoustic-tensor} are simple if $|\alpha|+|\beta|+|\gamma|\ll 1$.
\end{proposition}
\begin{proof}
Consider the perturbation $t(\xi)$ on the isotropic $\mu$-eigenspace:
\[ \pi_s(\xi)t(\xi)\pi_s(\xi) = \beta(J|\xi)^2 \pi_s(\xi) +(\beta \xi^2+\gamma(J|\xi)^2) \pi_J(\xi), \]
$\pi_J(\xi)$ the orthoprojector onto the line spanned by $\pi_s(\xi)J$.
By assumption, $\pi_s(\xi)t(\xi)\pi_s(\xi)$ has two distinct eigenvalues on $\image \pi_s(\xi)$.
It follows from \cite[Theorem \RN 2-5.4]{Kato76perturb} that, for small $\kappa\neq 0$,
the eigenvalues of $\ell_{\lambda,\mu}(\xi)+\kappa t(\xi)$ are pairwise distinct.
\end{proof}

\section{Spectral factorization and impedance}
\label{sect-spectral-factorization}

We decompose a quadratic self-adjoint matrix-valued polynomial
into a product first order factors which have disjoint spectra.
Real eigenvalues are classified according to their sign characteristic, \cite{GohbergLancRodman82matpolyn}.
Using ideas from \cite{MarkusMatsaev87oppencil} and simplifying assumptions,
we give complete proofs of the spectral factorization results needed in the main part of the paper,
without reference to the full sign characteristic.

Let $E$ be a finite-dimensional complex Hilbert space.
The scalar product $(x|y)$ is linear in the $x$ and conjugate linear in the $y$ variable.
We denote by $A^*$ the adjoint of linear map $A$.
Let $A(s)$ be a selfadjoint quadratic polynomial function with values in $\End(E)$:
\[ A(s)=A_0s^2+(A_1+A_1^*)s+A_2, \quad s\in\C;\]
$A(s)^*=A(\bar s)$, and $A_0$ is positive definite.
We are interested in factorisations
\begin{equation}
\label{sqp-factorisation}
A(s)=(s-Q^\sharp)A_0(s-Q).
\end{equation}
with right roots $Q\in \End(E)$.
The spectrum $\spec(A)$ is the set of $s\in\C$ for which $A(s)$ is singular.
The spectrum is finite and invariant under complex conjugation.
The non-zero elements of the null-space $\kernel A(s)$ are called eigenvectors associated with the eigenvalue $s$.
Given $\sigma\subset\C$, a factorization \eqref{sqp-factorisation} is said to be
$\sigma$-spectral if $\spec(Q)=\spec(A)\cap\sigma$ and $\spec(Q^\sharp)=\spec(A)\setminus\sigma$.
For a $\sigma$-spectral right root $Q$ of \eqref{sqp-factorisation},
$\kernel A(s)=\kernel(s-Q)$ holds if $s\in\sigma$.
Furthermore, the root polynomials at $\sigma$ of $A(s)$ and $s-Q$ are the same.
We associate the impedance $Z=-i(A_0Q+A_1)\in\End(E)$ with a factorization \eqref{sqp-factorisation}.

A straightforward computation shows that the solvency equation
\begin{equation}
\label{sqp-solvency}
A_0Q^2+(A_1+A_1^*)Q+A_2=0
\end{equation}
is equivalent to \eqref{sqp-factorisation}.
The left root $Q^\sharp$ is determined from $Q^\sharp A_0+A_0Q+A_1+A_1^*=0$.
Suppose $Q$ is a $\sigma$-spectral right root $Q$ of \eqref{sqp-factorisation},
and $\gamma$ a closed Jordan curve containing $\spec(Q^\sharp)$ in its exterior.
Cauchy's Theorem applied to $(s -Q)A(s)^{-1}$ gives
\begin{equation}
\label{sqp-Qr-intformula}
Q \oint_\gamma A(s)^{-1}\intd s
   = \oint_\gamma s A(s)^{-1}\intd s.
\end{equation}
This formula determines $Q$ on the range of $\oint_\gamma A(s)^{-1}\intd s$.

Associate with $A(s)$ the matrix
\begin{equation}
\label{factor-Stroh-matrix}
S = \begin{bmatrix} -A_0^{-1}A_1&A_0^{-1}\\ -A_2+A_1^*A_0^{-1}A_1& -A_1^*A_0^{-1}\end{bmatrix}\in\End(\EE).
\end{equation}
A simple calculation shows that
\begin{equation}
\label{factor-Stroh-linearization}
\begin{bmatrix}s A_0+A_1^*& I\\ -A_0& 0\end{bmatrix} (s-S)
   = \begin{bmatrix}A(s)& 0\\ -s A_0 -A_1 & I\end{bmatrix}.
\end{equation}
Put $J_1 = \begin{bmatrix}I&0\end{bmatrix}, J_2 = \begin{bmatrix} 0& I \end{bmatrix} \in\Hom(\EE,E)$.
Passing to inverses in \eqref{factor-Stroh-linearization}, we get
\begin{equation}
\label{sqp-Stroh-resolvent}
(s-S)^{-1} J_2^* = \begin{bmatrix} A(s)^{-1} \\ (s A_0+A_1) A(s)^{-1} \end{bmatrix}.
\end{equation}
In particular, $A(s)^{-1}=J_1(s-S)^{-1} J_2^*$.
Thus $\spec(A)=\spec(S)$, and the orders of the poles of $(s-S)^{-1}$ and of $A(s)^{-1}$ agree.
The linear polynomial $s-S$ is called a linearization of the matrix polynomial $A(s)$,
\cite{GohbergLancRodman82matpolyn}.

\begin{remark}
The special linearization \eqref{factor-Stroh-matrix} is taken from the Stroh formalism of elasticity,
see \cite{Stroh62} and \cite[Theorem~1.2]{Tanuma07Stroh}.
The second order ODE systems $A(D)u=0$ with symbol $A(s)$, $D=-i\intd/\intd r$,
is equivalent to the first order system $(D-S)w=0$, where $w(r)=[u(r),(A_0D+A_1)u(r)]^t$
is the joint field of displacements and tractions.
The components of $w(0)$ are related through the boundary impedance $Z$.
Suppose $Q$ is a $\sigma$-spectral right root of $A(s)$
such that $\sigma\subset\{\IM s\geq 0\}$ holds and such that the real eigenvalues are semisimple.
The solutions of the first order system $(D-Q)u=0$, that is $u(r)=e^{irQ}u(0)$,
are solutions of $A(D)u=0$ which are bounded on the positive half-line, $r\geq 0$.
\end{remark}

Applying the Dunford--Taylor functional calculus of the operator $S$, we obtain
\begin{equation}
\label{sqp-Stroh-funccalc}
(2\pi i)^{-1} \oint f(s) A(s)^{-1}\intd s = J_1 f(S) J_2^*
\end{equation}
if $f$ is holomorphic in a neighborhood of $\spec(A)$.
The contour of integration is a closed curve in the domain of $f$ which winds singly
around each point of $\spec(A)$.
In particular, for $\sigma\subset\spec(A)$, we have the projector $P_\sigma=1_\sigma(S)$
onto the sum of generalized eigenspaces corresponding to $\sigma$.
Here $1_\sigma=1$ in a neighborhood of $\sigma$ and $=0$ in a neighborhood of $\spec(A)\setminus\sigma$.

Denote the derivative of $A(s)$ by $A'(s)$.
\begin{lemma}
\label{lemma-semisimple-eigenval}
Let $s\in\spec(A)$.
The following conditions are equivalent:
\begin{compactenum}[(a)]
\item\label{item-simplicity-mu}
$A'(s)v+A(s)v_1\neq 0$ if $0\neq v\in\kernel A(s)$ and $v_1\in E$.
\item\label{item-simplicity-mu-R}
There exist $R,R_1\in\End(E)$ such that $RA(s)=0$ and $RA'(s)+R_1A(s)=I$.
\item\label{item-simplicity-mu-pole}
The pole of $A(s)^{-1}$ at $s$ is of order $1$.
\end{compactenum}
\end{lemma}
We denote by $R_s$ the residue at $s$ of the meromorphic map $A(s)^{-1}$.
\begin{proof}
Each of the conditions \eqref{item-simplicity-mu} and \eqref{item-simplicity-mu-R}
is equivalent to the map $A'(s)$ inducing a bijection from $\kernel A(s)$
onto the quotient $E/\image A(s)$.
Condition \eqref{item-simplicity-mu-pole} implies \eqref{item-simplicity-mu-R} with $R=R_s$.
Suppose the order of the pole is $>1$.
Then there exists a non-zero principal vector $w'$, that is,
we have $(s-S)w'+w=0$ and $(s-S)w=0$ for some $w\neq 0$.
The first component $v$ of the eigenvector $w$ is non-zero and satisfies $A(s)v=0$,
the second component equals $(s A_0+A_1)v$.
Using the first row of \eqref{factor-Stroh-linearization},
we find that the first component of $w'$, $v'$, satisfies
$A(s)v'+(2s A_0+A_1^* +A_1)v=0$.
This contradicts \eqref{item-simplicity-mu}.
\end{proof}
Suppose $s$ is semisimple, that means, the conditions of Lemma~\ref{lemma-semisimple-eigenval} hold.
Then $\image R_s=\kernel A(s)$ and $\image A(s)=\kernel R_s$.

Suppose $s\in\R\cap\spec(A)$ and $(A'(s)v|v)\neq 0$ if $0\neq v\in\kernel A(s)$.
Since $(A'(s)v|v)$ is real and $\kernel A(s)\setminus 0$ connected,
the quadratic form on $\kernel A(s)$ defined by $A'(s)$
is either positive definite or negative definite.
We say that $s$ is of \emph{positive type} (resp.\ \emph{negative type}) iff
$(A'(s)v|v)$ is positive (resp.\ negative) for all $v\in\kernel A(s)$, $v\neq 0$.

\begin{lemma} [{\cite[Lemma~2.1]{MarkusMatsaev87oppencil}}]
\label{sqp-residue-semidef}
Suppose $s\in\R\cap\spec(A)$ is of positive (resp.\ negative) type.
Then $s$ is semisimple.
The residue $R_s$ is positive (resp.\ negative) semidefinite.
\end{lemma}
\begin{proof}
Let $v\in \kernel A(s)$ and $v_1\in E$.
Thus
$\big(A'(s)v+A(s)v_1\big| v\big) = (A'(s)v|v)$.
By the assumption and \eqref{item-simplicity-mu}, $s$ is semisimple.
Observe that $R_s$ vanishes on the orthogonal complement of $\kernel A(s)$.
Condition \eqref{item-simplicity-mu-R} with $R=R_s$ implies $(A'(s)v|R_s v)=(v|v)$.
Hence the eigenvalues of $R_s|_{\kernel A(s)}$ are positive (resp.\ negative)
iff $s$ is of positive (resp.\ negative) type.
\end{proof}

The following result is a special case of \cite[Theorem~11.4]{GohbergLancRodman82matpolyn}.
\begin{proposition}
\label{prop-sqp-factorization}
Suppose each real eigenvalue of $A(s)$ is of positive or negative type.
Let $\sigma\subset\spec(A)$ be such that complex conjugation maps $\sigma\setminus\R$
bijectively onto $\sigma^\sharp\setminus\R$, where $\sigma^\sharp=\spec(A)\setminus\sigma$.
Assume that the eigenvalues in $\sigma\cap\R$ are all of the same type, opposite to the type
of the eigenvalues in $\sigma^\sharp\cap\R$.
Put $M=P_\sigma (\EE)$.
Let $\gamma$ be a positively oriented closed Jordan curve which contains $\sigma$
in its interior and $\sigma^\sharp$ in its exterior.
The following hold:
\begin{compactenum}[(i)]
\item
$J_1$ restricts to a linear isomorphism $J:M\to E$, and $J^{-1}=[I,iZ]^t$.
\item
$Q=JS|_M J^{-1}$ is the right root of a unique $\sigma$-spectral factorization \eqref{sqp-factorisation}.
\item
$\oint_\gamma A(s)^{-1}\intd s$ is non-singular.
\end{compactenum}
\end{proposition}
\begin{proof}
Write $P=P_\sigma$.
Put $C_j=J_1 PS^j P J_2^* = \oint_\gamma s^j A(s)^{-1}\intd s/2\pi i$.
Here we use \eqref{sqp-Stroh-funccalc}.
\begin{equation}
\label{eq-RExCjx-zero-implies-x-zero}
\RE( x|C_0 x) = \RE( x|C_1 x) = 0 \implies x=0
\end{equation}
holds for $x\in E$.
To see this, put $\bar\gamma(t)= \overline{\gamma(t)}$.
Since the polynomial $A$ is self-adjoint,
\[
C_j+C_j^* =(2\pi i)^{-1}\oint_{\gamma-\bar\gamma} s^j A(s)^{-1}\intd s.
\]
Note that $A(s)^{-1}=A_0^{-1}s^{-2}+\bigoh(|s|^{-3})$ as $s\to\infty$.
Applying Cauchy's Theorem, we get
\begin{align*}
C_0+C_0^* &= \sum\nolimits _{s\in\sigma\cap\R} R_s - \sum\nolimits _{s\in\sigma^\sharp\cap\R} R_s, \\
C_1+C_1^* &= A_0^{-1} + \sum\nolimits _{s\in\sigma\cap\R} s R_s
             - \sum\nolimits _{s\in\sigma^\sharp\cap\R} s R_s.
\end{align*}
The first equation implies
\[
2\RE( x|C_0 x)
  = \sum\nolimits _{s\in\sigma\cap\R} ( x| R_s x) -\sum\nolimits _{s\in\sigma^\sharp\cap\R} ( x| R_s x).
\]
Suppose $\RE( x|C_0 x) = 0$.
In view of Lemma~\ref{sqp-residue-semidef} and our assumptions,
the right-hand side is a sum of terms which are all of the same sign, hence each term is zero.
It follows that $2\RE( x|C_1 x) = ( x|A_0^{-1} x)$, which implies $x=0$ if $\RE(x|C_1x)=0$.

If $x$ is orthogonal to $J_1M=J_1P(\EE)$, then $(x|C_jx)=0$, hence $x=0$ by \eqref{eq-RExCjx-zero-implies-x-zero}.
Therefore, $E=J_1M$.
Similarly, replacing $\sigma$ by $\sigma^\sharp$, we get $E=J_1M^\sharp$ with $M^\sharp=P_{\sigma^\sharp}(\EE)$.
Note that $\EE=M\oplus M'$.
Thus $\dim M=\dim M^\sharp=\dim E$.
We infer that $J=J_1|M$ maps $M$ bijectively onto $E$.

Suppose $x\in E$ is orthogonal to the range of $J_2P^*$.
In view of $C_j^*=J_2 P^*(S^j)^* P^* J_1^*$, there holds $(x|C_j^*x)=0$ for all $j$.
Again using \eqref{eq-RExCjx-zero-implies-x-zero}, we infer $x=0$.
Hence $J_2P^*$ is surjective, thus $PJ_2^*$ injective and onto $M$.
Therefore $C_0=J_1PPJ_2^*$ is non-singular.

Observe that $Q=JS|_M J^{-1}\in \End(E)$ satisfies $C_j= Q^j C_0$.
Therefore,
\[
A_0Q^2C_0+(A_1+A_1^*)QC_0+ A_2C_0
  = (2\pi i)^{-1}\oint_\gamma A(s)A(s)^{-1}\intd s =0.
\]
This implies the solvency equation \eqref{sqp-solvency},
hence a factorization \eqref{sqp-factorisation}.

From \eqref{factor-Stroh-linearization} we derive
$\det(s-S) = \det(A_0)^{-1}\det(A(s))$.
Using \eqref{sqp-factorisation} this leads to
\[ \det(s-S) =\det(s-Q^\sharp)\det(s-Q). \]
Note that $Q$ is similar to $S|_M$, and $S$ is similar to $S|_M\oplus S_{M^\sharp}$.
Hence $\det(s-S|_{M^\sharp}) =\det(s-Q^\sharp)$,
implying that the factorization \eqref{sqp-factorisation} is $\sigma$-spectral.

In view of \eqref{sqp-Qr-intformula} and the invertibility of $C_0$,
it is clear that $Q$ is the only $\sigma$-spectral right root of $A(s)$.

Integrating \eqref{sqp-Stroh-resolvent} over $\gamma$, we obtain
$P_\sigma J_2^* = [I,A_0Q+A_1]^t C_0$.
Therefore, $[I,A_0Q+A_1]^t$ maps into $M$, and is a right inverse of $J$, hence the inverse of $J$.
\end{proof}

Suppose $\sigma\subset\spec(A)$ satisfies the assumptions in Proposition~\ref{prop-sqp-factorization}
and $\sigma\subset\{\IM s\geq 0\}$.
If $\spec(A)$ contains no real eigenvalues, then $\sigma$ is uniquely determined.
If there are real eigenvalues, then there are exactly two possibilities for $\sigma$:
The real eigenvalues in $\sigma$ are either of positive or of negative type.
We distinguish the cases by saying that $\sigma$ is of positive or negative type,
$\sigma=\sigma_+$ or $\sigma=\sigma_-$.
Denote by $Q$ the right spectral root of $A(s)$ determined by $\sigma$,
and by $Q_\pm$ the right spectral roots determined by $\sigma_\pm$.
The self-adjointness of the polynomial $A$ implies the formula $Q^{\sharp}_\pm=Q^{*}_\mp$,
a relation between left roots and adjoints of right roots.

The isomorphism $J:M\to E$ is a similarity between $S|_M$ and $Q$.
Put $M_c=P_{\sigma\setminus\R}M$ and $E_c=JM_c$.
The real eigenvalues $s$ of $S|_M$ are semisimple, and $J$ maps the eigenspaces onto $\kernel(Q-s)=\kernel A(s)$.
We have the following decomposition of $E$ into $Q$-invariant subspaces:
\begin{equation}
\label{eq-H-spec-decomp}
E = E_c \oplus E_r, \quad E_r = \oplus_{s\in\sigma\cap\R} \kernel A(s).
\end{equation}
Observe that $E_c$ equals the range of $\oint_\gamma A(s)\intd s$ whenever
$\gamma$ is a closed Jordan curve in the upper half-plane which contains $\sigma\setminus\R$ in its interior.
Applying \eqref{sqp-Qr-intformula}, we see that $Q|_{E_c}$ does not depend on $\sigma\cap\R$
and that $\spec(Q|_{E_c})\subset\{\IM s>0\}$.

Next we study the impedance $Z$.
Eigenvectors of $Q$ give explicit formulas:
\begin{equation*}
Qv=s v \implies iZv=(s A_0+A_1)v.
\end{equation*}
Below, we shall see that the restrictions of $Z$ to $E_c$ and $E_r$ behave rather differently.

Suppose $v_j\in E_c$.
Set $u_j(r)=e^{irQ}v_j$ for $r\geq 0$.
Then
\[ (Zv_1|v_2) = -\int_0^\infty D \big((A_0D+A_1)u_1(r)\big|u_2(r)\big)\intd r, \]
$D=-i\intd/\intd r$.
Using the differential equation $A(D)u_1=0$, we find that
\begin{equation}
\label{eq-Z-int-over-r}
(Zv_1|v_2) = \int_0^\infty (A_0Du_1|Du_2)+(A_1u_1|Du_2)+(Du_1|A_1u_2)+(A_2u_1|u_2)\intd r.
\end{equation}
As $A_0$ and $A_2$ are self-adjoint,
this implies that $(Zv_2|v_1)$ is the complex conjugate of $(Zv_1|v_2)$.
Hence $(Z^*v_1|v_2)=(Zv_1|v_2)$.
Thus we have shown:
\begin{equation}
\label{eq-Z-selfadj-Ec}
Z|_{E_c}=Z^*|_{E_c}
\end{equation}
It follows that $\IM(u|Zu)=0$ holds for $u\in E_c$.
In fact, $\kernel Z\subset E_c$ by the following.

\begin{proposition}
\label{prop-Z-nonsing-if-sigma-real}
Let $u\in E$.
Denote by $u_s$ the $\kernel A(s)$-component of $u$ with respect to \eqref{eq-H-spec-decomp}.
There holds
\begin{equation}
\label{eq-IM-u-Zu}
2\IM(u|Z u) = \sum_{s\in\sigma \cap\R} (A'(s)u_s|u_s).
\end{equation}
If $\sigma$ has positive (resp.\ negative) type, then
$\IM(u|Z u)\geq 0$ (resp.\ $\IM(u|Z u)\leq 0$) holds for $u\in E$,
with strict inequality if $u\not\in E_c$.
\end{proposition}
\begin{proof}
We use the following indefinite inner product on $\EE$:
\[
[w_1|w_2] =(w_1|Kw_2)_{\EE}=(x_1|y_2) + (y_1|x_2),
\quad
w_j = \begin{bmatrix} x_j\\ y_j \end{bmatrix}\in \EE, \quad
K = \begin{bmatrix} 0&I\\ I&0 \end{bmatrix}.
\]
Directly from the definition \eqref{factor-Stroh-matrix} we see that
the map $S$ is self-adjoint with respect to the inner product, that is,
$KS$ is selfadjoint with respect to the standard scalar product of $\EE$.
Therefore, $[(S-s)^r w|v]= [w|(S-\bar s)^rv]$.
It follows that
\begin{equation*}
[w|v]=0 \quad
\text{if $(S-s)^r w=0=(S-s)v$ and $\bar s\neq s$.}
\end{equation*}
Hence
$[w|w]=[\tilde w|\tilde w]+\sum\nolimits _{s\in\sigma\cap\R} [w_s|w_s]$
if $w=\tilde w+\sum_s w_s$, $\tilde w\in M_c$, $w_s\in \kernel(S|_M-s)$.
We apply the formula with $w=(u,iZ u)$.
Note that $[w|w]=2\IM(u|Z u)$.
By \eqref{eq-Z-selfadj-Ec}, $[w|w]=0$ if $u\in E_c$.
We find that
\[ \IM(u|Z  u) = \sum\nolimits_{s\in\sigma\cap\R} \IM(u_s|Z  u_s). \]
Suppose $u\in\kernel A(s)$, $s\in\sigma\cap\R$.
Then
\[ (A'(s)u|u)= ((2s A_0+A_1+A_1^*)u|u) = (u|iZ u)+(iZ u|u)=2\IM(u|Z u). \]
Summarizing, we have shown \eqref{eq-IM-u-Zu}.
To complete the proof,
recall that $\pm (A'(s)u|u)>0$ holds for $0\neq u\in\kernel A(s)$ if $s\in\R\cap\sigma_\pm$.
\end{proof}

Finally, we present and reprove results which stem from and underlie the Barnett--Lothe theory
of (Rayleigh) surface waves in anisotropic elastic media,
\cite{LotheBarnett85surfwaveimped}.
We assume that $A(s)$ has no real eigenvalues.
So we consider the $\sigma$-spectral factorization for $\sigma=\spec(A)\cap\{\IM s>0\}$.
Therefore $E=E_c$, and $Z$ is self-adjoint by \eqref{eq-Z-selfadj-Ec}.
By \eqref{sqp-Qr-intformula} and the definition of $Z$ we have
\begin{equation*}
iZ \oint_\gamma A(s)^{-1}\intd s
    = \oint_\gamma (s A_0 +A_1)A(s)^{-1}\intd s.
\end{equation*}
For $R>0$ sufficiently large, we take the closed curve $\gamma$ as the sum of the interval $[-R,R]$
and the semicircle $Re^{it}$, $0\leq t\leq \pi$.
Letting $R\to\infty$ and using
\(s A_0 A(s)^{-1}=s^{-1}\Id +\bigoh(s^{-2}),\)
we get the Barnett--Lothe integral formula:
\begin{equation}
\label{eq-Z-int-BL}
iZ \int_{-\infty}^\infty A(s)^{-1}\intd s
    = i\pi \Id +\lim_{R\to\infty}\int_{-R}^R (s A_0 +A_1)A(s)^{-1}\intd s.
\end{equation}
Note that $A(s)$ is positive definite for real $s$.

Real and imaginary parts of linear operators on $E$ are defined if the Hilbert space $E$ carries a real structure.
This means that there exists a complex conjugation involution on $E$,
and $E$ is the complexification of the real subspace on which the conjugation restricts to the identity.

\begin{proposition}
\label{prop-Re-Z-posdef}
Suppose the Hilbert space $E$ carries a real structure, and that the operators $A_j$ are real.
Suppose $\spec(A)\cap\R=\emptyset$ and $\sigma=\spec(A)\cap\{\IM s > 0\}$.
The real part $\RE Z$ of $Z$ is positive definite.
In addition, suppose $\dim E\geq 3$.
There exists a $2$-dimensional $Z$-invariant subspace of $E$ on which $Z$ is positive definite.
\end{proposition}
\begin{proof}
From \eqref{eq-Z-int-BL} we infer that $\RE Z=\pi\big(\int_{-\infty}^\infty A(s)^{-1}\intd s\big)^{-1}$,
which proves the positive definiteness of $\RE Z$.
Suppose $\dim E\geq 3$ holds, and $Z$ has at most one positive eigenvalue.
Then there exists $0\neq w\in E$ such that $(v|Zv)\leq 0$ whenever $(v|w)=0$.
Let $v\neq 0$ be a real vector which is orthogonal to both $\RE w$ and $\IM w$.
Then $(v|(\RE Z)v)=(v|Zv)\leq 0$ which contradicts the positive definiteness of $\RE Z$.
\end{proof}

\begin{proposition}
\label{prop-Zdot-posdef}
Suppose $\spec(A)\cap\R=\emptyset$ and $\sigma=\spec(A)\cap\{\IM s > 0\}$.
Suppose $A_2$ depends smoothly on a real parameter, whereas $A_0$ and $A_1$ do not depend on that parameter.
Denoting the derivative with respect to the parameter by a dot,
$\dot Z$ is positive or negative definite if $\dot A_2$ is.
\end{proposition}
\begin{proof}
Insertion of $Q=A_0^{-1}(iZ-A_1)$ into the solvency equation \eqref{sqp-solvency} gives
\begin{equation*}
(iZ+A_1^*)A_0^{-1}(iZ-A_1) +A_2 =0,
\end{equation*}
which is an equation of Ricatti type, \cite{MielkeFu04surfwavespeed}.
Differentiating, we obtain the Lyapunov equation
$i(\dot Z Q-Q^*\dot Z)  = - \dot A_2$.
Since the spectra of $Q$ and $Q^*$ are disjoint, the equation is uniquely solvable.
Using
$\int_0^\infty \frac{\intd}{\intd r} e^{irQ^*}Y e^{riQ}\intd r = -Y$,
we obtain $\dot Z=\int_0^\infty e^{irQ^*}\dot A_2 e^{riQ}\intd r$.
\end{proof}

\def\cprime{$'$}
\providecommand{\MRev}[1]{%
  \href{http://www.ams.org/mathscinet-getitem?mr=#1}{MR #1}
}
\providecommand{\Zbl}[1]{%
  \href{https://zbmath.org/?q=an:#1}{Zbl #1}
}
\providecommand{\arXiv}[2]{%
  \href{http://arxiv.org/abs/#1}{arXiv:#1 [#2]}
}
\providecommand{\href}[2]{#2}

\end{document}